\documentclass[sn-mathphys,Numbered,icol]{sn-jnl}

\usepackage{graphicx}
\graphicspath{{./images}}
\usepackage{multirow}
\usepackage{amsmath,amssymb,amsfonts}
\usepackage{amsthm}
\usepackage{mathrsfs}
\usepackage[title]{appendix}
\usepackage{xcolor}
\usepackage{textcomp}
\usepackage{manyfoot}
\usepackage{booktabs}
\usepackage{algorithm}
\usepackage{algorithmicx}
\usepackage{algpseudocode}
\usepackage{listings}
\usepackage{cleveref}
\usepackage{ulem}

\usepackage{tkz-euclide}
\usepackage{tikz,pgf,pgfplots}
\usepackage{standalone}

\theoremstyle{thmstyleone}

\theoremstyle{thmstyletwo}

\theoremstyle{thmstylethree}

\begin{document}

\title[Computational homogenization for {{porous materials}}]{Computational homogenization for {{aerogel-like polydisperse open-porous materials}} using neural network--based surrogate models on the microscale}

\author*[1,2]{\fnm{Axel} \sur{Klawonn}}\email{axel.klawonn@uni-koeln.de}
\equalcont{List of authors is in alphabetical order.}

\author*[1,2]{\fnm{Martin} \sur{Lanser}}\email{martin.lanser@uni-koeln.de}
\equalcont{List of authors is in alphabetical order.}

\author*[1]{\fnm{Lucas} \sur{Mager}}\email{lucas.mager@uni-koeln.de}
\equalcont{List of authors is in alphabetical order.}

\author*[3,4]{\fnm{Ameya} \sur{Rege}}\email{ameya.rege@dlr.de}
\equalcont{List of authors is in alphabetical order.}

\affil[1]{\orgdiv{Department of Mathematics and Computer Science}, \orgname{University of Cologne}, \orgaddress{\street{Weyertal 86-90}, \city{Cologne}, \postcode{50931}, \country{Germany}}}

\affil[2]{\orgdiv{Center for Data and Simulation Science}, \orgname{University of Cologne}, \orgaddress{\street{Albertus-Magnus-Platz}, \city{Cologne}, \postcode{50923}, \country{Germany}}}

\affil[3]{\orgdiv{Institute of Materials Research}, \orgname{German Aerospace Center}, \orgaddress{\street{Linder H\"ohe}, \city{Cologne}, \postcode{51147}, \country{Germany}}}

\affil[4]{\orgdiv{{School of Computer Science and Mathematics}}, \orgname{Keele University}, \orgaddress{\street{The Covert}, \city{Staffordshire}, \postcode{ST5 5BG}, \country{United Kingdom}}}

\abstract{The morphology of nanostructured materials exhibiting a polydisperse porous space, such as aerogels, is very open porous and fine grained.
Therefore, a simulation of the deformation of a large aerogel structure resolving the nanostructure would be extremely expensive. Thus, multi-scale or homogenization approaches have to be considered. Here, a computational scale bridging approach based on the FE$^2$ method is suggested, where the macroscopic scale is discretized using finite elements while the microstructure of the open-porous material is resolved as a network of Euler-Bernoulli beams. Here, the beam frame based RVEs (representative volume elements) have pores whose size distribution follows the measured values for a specific material. This is a well-known approach to model aerogel structures. For the computational homogenization, an approach to average the first Piola-Kirchhoff stresses in a beam frame by neglecting rotational moments is suggested. To further overcome the computationally most expensive part in the homogenization method, that is, solving the  RVEs and averaging their stress fields, a surrogate model is introduced based on neural networks. The networks input is the localized deformation gradient on the macroscopic scale and its output is the averaged stress for the specific material. It is trained on data generated by the beam frame based approach. The effiency and robustness of both homogenization approaches is shown numerically, the approximation properties of the surrogate model is verified for different macroscopic problems and discretizations. Different (Quasi-)Newton solvers are considered on the macroscopic scale and compared with respect to their convergence properties.}


\maketitle

\section{Introduction}\label{sec1}

{{The mechanics of open-porous materials form a fascinating branch of materials science, unveiling a complex interplay of structural intricacies and mechanical behaviors that have profound implications across multiple disciplines. These materials, characterized by their interconnected network of voids or pores, exhibit unique mechanical properties driven by their open-porous architecture \cite{gibson-ashby-1999}. The mechanical performance is influenced by factors such as pore size, shape, distribution, and the material composition itself. Understanding the deformation mechanisms and stress distribution within these porous structures is crucial for optimizing their performance in applications ranging from lightweight structural components to advanced filtration systems. A special class of such materials is the territory of nanoporous materials, such as nanofoams and aerogels. Aerogels exhibit pore sizes ranging from as low as less than two nanometers to as high as over hundreds of nanometers within the same material \cite{fricke-emmerling-1992}. Such a polydisperse nature significantly affects the thermal conductivity and sound absorption characteristics, but also their overall mechanical performance \cite{rege-handbook-2023}. Formulating an intricate understanding of the macroscopic mechanics of such polydisperse open-porous materials demands the theoretical understanding of the underlying deformation mechanisms. While state-of-the-art experimental tools help to dive deep into the microstructure of the material, computational modeling has proven to establish concrete structure-property relations providing, generally, a physics-informed explanation.

Diverse studies have reported on the computational modeling of open-porous materials that exhibit pores on the macroscale, however, very few studies have investigated the computational description of nanostructured open-porous materials}. Here, we focus and limit our search to physics-informed or micromechanical models. Existing literature uses either constitutive modeling or computational homogenization as a means to describe the mechanical behavior of open-porous materials. Such materials are made up of a network of three-dimensionally interconnected struts or pore-walls. One of the first reports on the mathematical description of the deformation in such materials was presented by Gent and Thomas \cite{gent-thomas-1959}, where they described the network of elastic foams to be formed of dumbbell-shaped elements called threads and dead volumes. Only extension of the threads, namely axial deformation in the network was considered. A much more robust theory was presented later by Gibson and Ashby, who modeled these struts as Euler-Bernoulli beams and provided a description of their bulk mechanical behavior \cite{gibson-prsa-1982, gibson-ashby-prsa-1982}. While there were extensions of these models presented, a notable development was presented by Warren and Kraynik \cite{warren-kraynik-1997} who modeled a perfectly ordered foam, whose unit cell was inspired by a regular tetrakaidecahedron referred to as the open-foam Kelvin cell. This accounted for a more complex geometrical representation of a unit cell as compared to the previously considered square or hexagonal cells in two dimensions, or cubic cells in three dimensions. This model was inspired by the work of Dement'ev and Tarakanov \cite{dementev-tarakanov-1970}, who studied plastic foams. While the Kelvin-cell model in \cite{warren-kraynik-1997} investigated only the elastic properties, the approach was extended to describe large deformations and crushing of foams, particularly metal foams using the finite-element-method \cite{jang-ijss-2010, gaitanaros-ijss-2012}. Interestingly, the model by Dement'ev and Tarakanov could show very good validations under large deformations with elastic foams. On the similar grounds of the strain energy approach, first proposed in \cite{dementev-tarakanov-1970}, Rege et al. \cite{rege-pre-2021} proposed a generalized micromechanical constitutive model, based on the representation of cell walls as beams, to describe the macroscopic stress-strain response of aerogel-like open-porous materials. The model was extended to capture the densification in a later study and showed very good validation with biopolymer aerogels \cite{rege-mat-2021}. These above-mentioned models, both constitutive or homogenization ones, base the models on a representative unit cell, strictly speaking one that represents 'a pore' and the surrounding pore walls. This is however insufficient to describe more polydisperse open-porous materials. To this end, the construction of a representative volume element (RVE) is necessary that accounts for the complete pore space and the diversity in structural features of the material in consideration. In such a case, the Voronoi-tessellation-approach is widely used. In such a representation, the network is made up of Voronoi cells, and the cell walls, representing struts, are modeled again as beams. Analysis of such a modeling approach in application to open-porous materials has been reported by several authors \cite{zhu-jmps-2001, zhu-ijss-2006, alsayednoor-harrison-mom-2013, alsayednoor-harrison-mom-2016, DLR_biopoly_aerogel2019}. These studies have however mostly dealt with two-dimensional problems. For three-dimensional modeling, most reports have used computer tomography images and reconstructed them for structure-property analysis \cite{Jang-ijss-2008, fu-jncs-2011, Skibinski-2017, foray-gels-2022}. While this provides a realistic picture of the microstructure, this approach becomes increasingly challenging as the pore sizes become smaller, particularly, when one is dealing with those below 50 nm, which is the case, for example in aerogels. Chandrasekaran et al. \cite{DLR_biopoly_aerogel2021} proposed a radical Voronoi method, wherein the pore space was represented by a random closed pack of polydisperse spheres and the Laguerre-Voronoi-tessellations were generated over the spheres. A heuristic analysis of the microstructural parameters was subsequently analyzed by Aney and Rege \cite{aney-rege-2023}. However, this approach of analyzing the RVE only helps understand the bulk behavior under some classical loading scenarios and does not account for geometrical effects that may arise due to the micro-macro transition. It is here where a multiscale homogenization approach is demanded and lacking in the literature for such highly open-porous polydisperse materials. This becomes increasingly important for large heterogeneous RVEs} because resolving their structure on the macroscopic scale using, for example, a beam model results in large linear systems which cannot be solved without enormous computational effort. To overcome this limitation, we suggest to use a computational homogenization approach based on the well-known FE$^2$ method~\cite{MieheBound,SMIT_Hom_1998,HacklSchroeder,Kouz_Hom_2001,FEYEL_FE2_1999}.

In the FE$^2$ approach, the microstructure of the considered material is only resolved locally on the second scale, usually named microscopic scale. In this article, instead of using finite element discretizations on both levels, we modify the FE$^2$ method and use beam frame models to resolve the aerogel structure on the mesoscale. This is a natural choice, since the use of beams to model aerogel structures on the mesoscopic scale is a state-of-the-art approach and already used and suggested in~\cite{DLR_biopoly_aerogel2016, DLR_biopoly_aerogel2018, DLR_biopoly_aerogel2021, abdusalamov-jpcb-2021}.
 Also in the field of lattice structures the analysis of microscopic beam networks is an important research topic~\cite{ma15020605}. Multiscale methods ~\cite{VIGLIOTTI201444, VIGLIOTTI201257,VIGLIOTTI2018231} as well as beam theory-based~\cite{Wang2001,Wang2004,DESHPANDE20011747,MASTERS1996403} and machine learning~\cite{KOEPPE2018147} approaches are used for the analysis of these materials. Honeycomb lattice structures are of special interest for this research since it is a quite common structural element in nature and engineering. However, methods for lattice structures usually cannot be applied in the case of aerogels due to their unstructured arrangement of the fibrils. In the recent years data-driven computing ~\cite{KIRCHDOERFER201681,Conti2018} has increased in popularity and has also found application in the simulation of porous materials~\cite{KORZENIOWSKI2022115487} including structure-property predictions of aerogels \cite{abdusalamov-sm-2021, pandit-srep-2024}.

In general, the FE$^2$ method is computationally demanding and to increase the computational efficiency and reduce the time to solution, parallel implementations have been used a lot in the past, cf.~\cite{fepar1,fepar2,fepar3}. An alternative to the extensive use of parallel computing is to drastically reduce the computational complexity of the FE$^2$ approach by training and exploiting machine learning based surrogate models on the microscale. Combinations of surrogate models based on deep neural networks (DNNs) and the FE$^2$ method are often denoted by DNN-FE$^2$ approaches and for different finite element applications, such methods have already been successfully studied in recent years. They have been proven to be robust and computationally very efficient in different contexts. Some of these approaches can be found in~\cite{surr2,surr3,surr4,surr5,surr6,surr7,surr8,surr9,surr10,surr11,wang2024}. Based on these existing DNN-FE$^2$ methods, we suggest a neural network (NN) based model which is applied locally in all integration points of the macroscopic finite element discretization to predict the average stress depending on the local macroscopic deformation. Therefore, the DNN replaces the solution of local beam frame problems and the averaging of stresses in the usual homogenization approach. In contrast to, e.g., ~\cite{surr2,surr3,surr4,surr5,surr6,surr7,surr8,surr9,surr10,surr11,wang2024,surr12}, we consider porous aerogel-like structures on the microscale and we use beam frame models to create the data to train the DNN surrogate model. 

To the best of our knowledge, in the present work, we extend the  large family of already existing homogenization approaches with an FE$^2$-related method based on beam frame models on the microscale and, for the first time, train deep learning-based surrogate models for this method. We also, for the first time, extensively discuss the use of Broyden–Fletcher–Goldfarb–Shanno (BFGS) \cite{broyden1970}, which is a good alternative if the consistent tangent is not available or expensive to use.

The remainder of the article is organized as follows. In section 2 we describe the computational homogenization approach using beam frame models on the microscale. We replace the beam frame models by a faster surrogate model based on NNs in section 3.  Finally, in section 4, we show numerical results and discuss the convergence of the different methods and the quality and accuracy of the results.

\section{Computational homogenization methods}\label{sec:methods}

Performing macroscopic simulations resolving the microscopic morphology of aerogels and similar porous materials on a single scale is computationally extremely expensive. Therefore, we consider a computational scale bridging approach where the morphological representation of the  open-porous material in consideration is resolved in localized representative volume elements (RVEs). Our approach is based on the well-known FE$^2$ method~\cite{MieheBound,SMIT_Hom_1998,HacklSchroeder,Kouz_Hom_2001,FEYEL_FE2_1999}. This method has been developed for the simulation of micro-heterogeneous solid materials as, for example, dual- or multi-phase steels. Usually, the macroscopic problem is discretized with comparably large finite elements while the finite element discretization of the RVEs resolves the microscopic structure. In our deviation from the FE$^2$ method, we also use a finite element discretization on the macroscopic scale but a beam frame model on the microscopic scale which suits the nanostructure of the aerogel. Following the FE$^2$ framework, in each Gaussian integration point of the macroscopic problem, one RVE of the aerogel structure discretized with beams is attached. Then, the RVE is deformed with respect to the local macroscopic deformation in the corresponding integration point and delivers, as a response, the macroscopic stress in that point by averaging over the stresses in the beam elements of the RVE. More details on the modeling of both scales are given in the following sections.

\subsection{Modeling the nanostructure of the RVE}

\begin{figure}
	\centering
	\includegraphics[width=5cm]{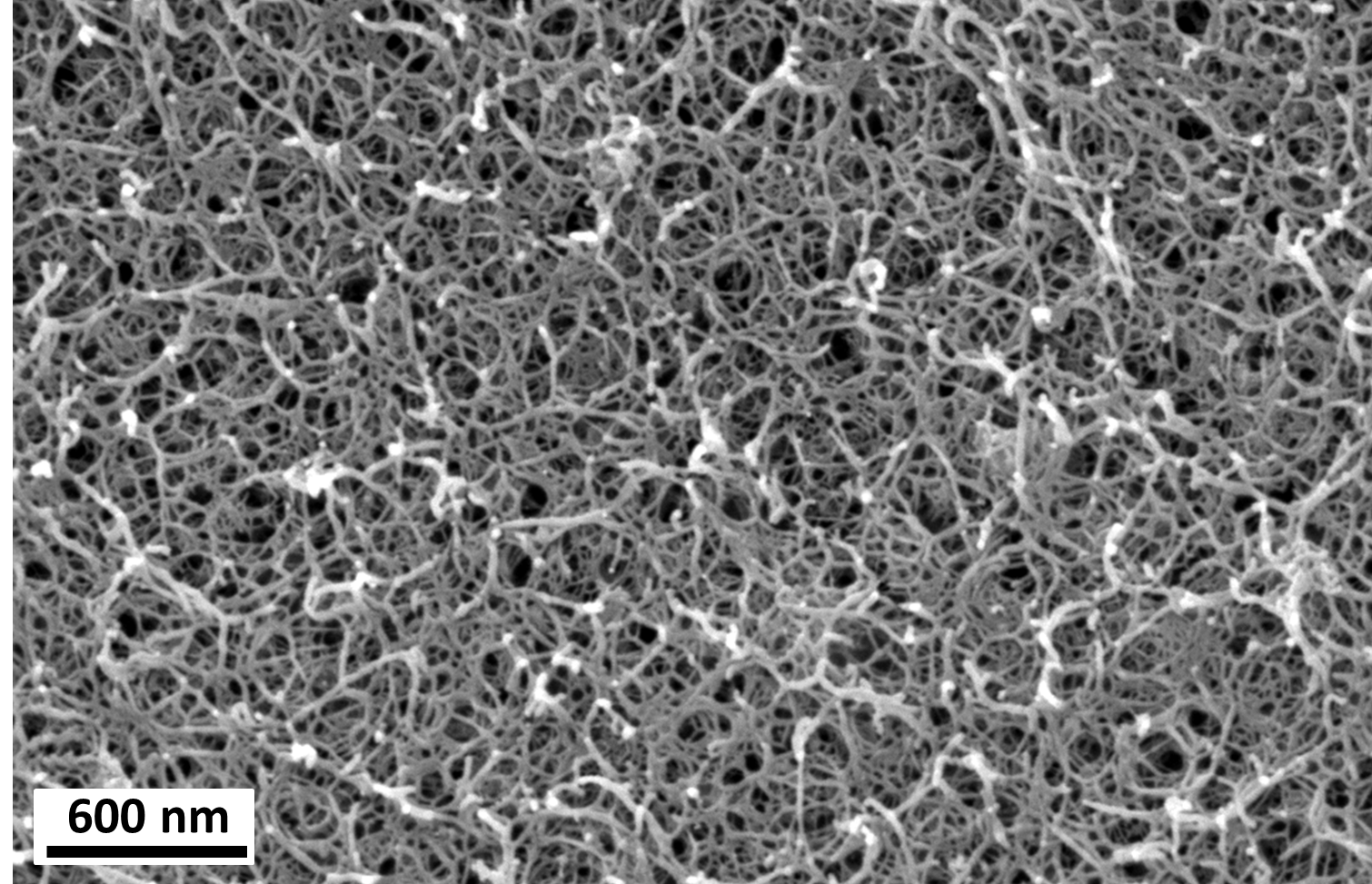}
	\caption{An SEM image of a biopolymer aerogel. The figure is adapted from Rege et al. \cite{DLR_biopoly_aerogel2018}.}
	\label{fig:aerogel_sem}
\end{figure}

The nanostructure of a typical open-porous material can be seen in \Cref{fig:aerogel_sem} which shows the scanning electron microscope (SEM) image of a carrageenan aerogel. For a two-scale homogenization approach, an RVE that is representative of the material's morphology, namely one that represents the entire pore space of the material in consideration, is desired. The RVE which represents such a nanostructure is generated using a method presented in \cite{DLR_biopoly_aerogel2019,DLR_biopoly_aerogel2021} for the case of biopolymer aerogels. The pore size distribution and the pore-wall diameter are required as an input for the generation algorithm. These parameters are obtained by experimental analysis such as SEM images and the Barrett-Joyner-Halenda (BJH) analysis of the physisorption isotherms and may vary depending on the exact procedure of the synthesis of the aerogel. A detailed description of the experimental analysis of the material can be found in \cite{DLR_biopoly_aerogel2018}.

\begin{figure}
\begin{minipage}{0.49\textwidth}
		\includegraphics[width=\textwidth]{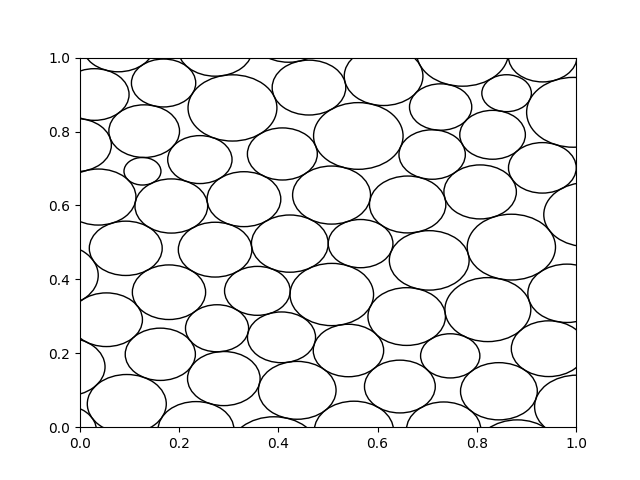}
\end{minipage}
\begin{minipage}{0.49\textwidth}
		\includegraphics[width=\textwidth]{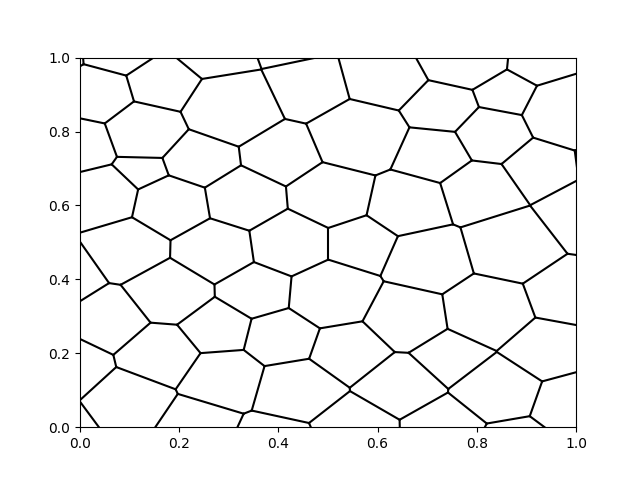}
\end{minipage}
\caption{Circle packing (left) and Voronoi tessellation (right) for the generation of the two-dimensional RVE.}
\label{fig:modelling}
\end{figure}

The first step of the RVE generation method is to use a sphere packing algorithm where the sphere volume distribution coalesces with the experimentally measured pore volume distribution of the given open-porous material. Afterwards, the sphere centers and the corresponding diameters are used in a Laguerre-Voronoi tessellation and the interfaces of the resulting Voronoi cells finally build the structure of the RVE. The struts on the interfaces represent the pore walls of the open-porous material. For a two-dimensional RVE the modelling is illustrated in \Cref{fig:modelling}. For more details, see also~\cite{DLR_biopoly_aerogel2019,DLR_biopoly_aerogel2021}.

\begin{figure}
	\centering
\includegraphics[width=0.35\textwidth]{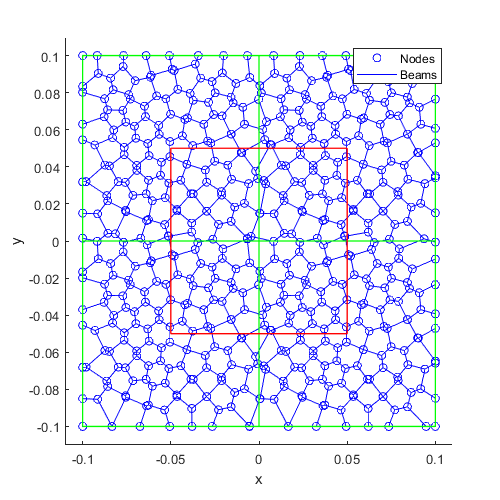}
\caption{Voronoi tessellation that results from copied sphere centers. Red square in the middle has periodic boundaries.}
\label{fig:cutout}
\end{figure}

To ensure periodic boundaries in each dimension it is possible to copy the spheres obtained from the sphere packing algorithm in each dimension and thus expand the domain. Cutting out the center of the Laguerre-Voronoi tessellation resulting from the increased number of sphere centers yields a periodic domain of the same size as the original domain and with continuous pore walls on each of the periodic boundaries. A visualization of this approach is presented in \Cref{fig:cutout} where the green squares mark the copied domains and the red square highlights the periodic boundaries of the center domain. Since the FE$^2$ method requires periodic boundary conditions on the microscopic scale we use this approach for the generation of the RVEs.

\subsection{Microscopic Problem}
On the microscopic level, the deformation of the pore walls in the RVEs is computed in each step of the overall homogenization algorithm and the resulting stress tensors are averaged afterwards. In the RVE, all the pore walls are modelled as Euler-Bernoulli beams and in the literature the resulting model is often related to as beam frame model~\cite{Oechsner2023}.

We consider an RVE consisting of $n_B$ beam elements and $n_N$ nodes which are all connected by the beams. The computation of the deformation and rotation of each of the beam vertices on the microscopic scale requires solving a linear system of equations $$K \cdot u=F.$$ Here, the stiffness matrix $K$ is assembled from beam stiffness matrices $K_e$ for each beam element $e$ \cite{kassimali2011}. We consider a single beam element $e=\left[v_i,v_j\right]$ with a circular cross section and with starting point and endpoint $$v_i = \left(\begin{array}{c}
	x_i \\
	y_i \\
	z_i \\
\end{array}\right) {\; \rm and \;} v_j = \left(\begin{array}{c}
	x_j \\
	y_j \\
	z_j \\
\end{array}\right).$$ Then, the beam stiffness matrix is given by

$$K_e = \mathbf{T}^T\cdot \frac{1}{L_e} \left(\begin{array}{llll}
	\mathbf{D} & \mathbf{C}^T & -\mathbf{D} & \mathbf{C}^T	\\
	\mathbf{C} & \mathbf{R_1} & \mathbf{C}^T & \mathbf{R_2}	\\
	-\mathbf{D} & \mathbf{C} & \mathbf{D} & \mathbf{C}	\\
	\mathbf{C} & \mathbf{R_2} & \mathbf{C}^T & \mathbf{R_1}
\end{array}\right)\cdot \mathbf{T},$$

\noindent where the submatrices $\mathbf{D}$, $\mathbf{R_1}$, and $\mathbf{R_2}$ are related to the deformation and the rotation of the beam element. The matrix $\mathbf{C}$ refers to the coupling of these variables. For the three-dimensional case the submatrices are defined as

\begin{align*}
	\mathbf{D} &= \left(\begin{array}{ccc}
		E\cdot A & 0 & 0 \\
		0 & 12 \frac{E\cdot I}{L_e^2} & 0    \\
		0 & 0 & 12 \frac{E\cdot I}{L_e^2}
	\end{array}\right),\\
	\mathbf{R}_1 &= \left(\begin{array}{ccc}
	\frac{E\cdot I^2}{2\cdot(1+\nu)} &  0 & 0  \\
	0 & 4 E\cdot I & 0    \\
	0 & 0 & 4 E\cdot I
\end{array}\right),	\\
\mathbf{R}_2 &= \left(\begin{array}{ccc}
	-\frac{E\cdot I^2}{2\cdot(1+\nu)} &  0 & 0  \\
	0 & 2 E\cdot I & 0    \\
	0 & 0 & 2 E\cdot I
\end{array}\right),\\ 
	{\rm and\;}  \mathbf{C} &= \left(\begin{array}{ccc}
	0 & 0 & 0 \\
	0 & 0 & -6\frac{E\cdot I}{L_e}    \\
	0 & 6\frac{E\cdot I}{L_e} & 0
\end{array}\right).
\end{align*}

Finally, $\mathbf{T}$ is a transformation matrix which depends on the orientation of the beam in space and the parameters $E$, $I$, $A$, and $L_e$ refer to Young's modulus, the second moment of area, the area of the cross section, and the length of the beam. For a detailed description of the beam frame model and the Euler-Bernoulli beam theory, we refer to~\cite{kassimali2011}.

The systems right hand side $F=\left(\mathbf{Q}_1,\mathbf{M}_1,\mathbf{Q}_2,\mathbf{M}_2\dots,\mathbf{Q}_{n_N},\mathbf{M}_{n_N}\right)^T$ consists of forces $\mathbf{Q}_i=\left(Q_i^x,Q_i^y,Q_i^z\right)^T$ and bending moments $\mathbf{M}_i=\left(M_i^x,M_i^y,M_i^z\right)^T$ at each node. Solving the equation yields the solution vector $u=\left(\mathbf{u}_1,\boldsymbol{\theta}_1,\mathbf{u}_2,\boldsymbol{\theta}_2,\dots,\mathbf{u}_{n_N},\boldsymbol{\theta}_{n_N}\right)^T$, where $\mathbf{u}_i=\left(u_i^x,u_i^y,u_i^z\right)^T$ refers to the deformation of each node and $\boldsymbol{\theta}_i=\left(\theta_i^x,\theta_i^y,\theta_i^z\right)^T$ is the vector of the corresponding rotations in each direction. By applying the beam frame model as described, it is assumed that the deformation within one beam element can be described by a cubic polynomial depending on its distance from the starting point of the beam. Let us introduce the variable $\xi\in[0,1]$ as a local measurement for each beam for the normalized distance from each of the ends. For the beam element $e=[v_i,v_j]$ the value $\xi=0$ for example refers to the global coordinates $v_i = \left(
	x_i,
	y_i,
	z_i \right)^T$ and the value $\xi=1$ refers to the other node $v_j = \left(
	x_j,
	y_j,
	z_j \right)^T$. With this local variable it is possible to describe the  deformation of each beam over its length in terms of $\xi$ as $u^{x_k}(x,y,z)=\Tilde{u}^{x_k}(\xi)$ for each $x_k=(x,y,z)$, where $\Tilde{u}$ is a polynomial of third degree.

In general, the deformation of each node $v_i$, $i=1,\dots,n_N$, in the RVE can be split into two parts

$$\mathbf{u}_i=\mathbf{\tilde{u}}_i+\mathbf{\bar{u}}_i.$$

Here, $\mathbf{\bar{u}}_i$ is defined by the macroscopic deformation $\overline{F}$ in the corresponding integration point using the relation $\mathbf{\bar{u}}_i = \overline{F}\cdot{v_i}$ while $\mathbf{\tilde{u}}_i$ is the microscopic fluctuation field, which is finally computed by solving the microscopic problem.

We apply periodic boundary conditions to the fluctuation and the rotations, that is,
$$ \mathbf{\tilde{u}}^+ = \mathbf{\tilde{u}}^-\quad \text{and } \boldsymbol{\theta}^+ = \boldsymbol{\theta}^-$$
for each pair of periodic nodes $v^+$ and $v^-$ at opposing faces of the RVE. To obtain a regular matrix $K$, in each corner of the RVE, the fluctuations are fixed with zero Dirichlet boundary conditions.

\subsection{Macroscopic Problem}

We formulate the macroscopic problem based on the weak formulation of the momentum balance equation which is given by
\begin{align*}
	\int_{\overline{B}_0} \delta\bar{x}\left({\rm Div}_{\bar{x}}\overline{P}(\overline{F})-\Bar{f}\right)d\Bar{x}=0
\end{align*}
for a test function $\delta\Bar{x}$. Without the consideration of any volume force the equation is reduced to

\begin{equation}
	\int_{\overline{B}_0} \delta\bar{x}\,{\rm Div}_{\bar{x}}\overline{P}(\overline{F})d\Bar{x}=0.
	\label{macro_problem}
\end{equation}

The principle behind the FE$^2$ method is that the relation between the deformation gradient $\overline{F}(\bar{x})$ and the macroscopic Piola-Kirchhoff stress $\overline{P}(\overline{F}(\bar{x}))$ in a certain point $\bar{x}$ is not modelled by a material law but by the results of microscopic simulations incorporating the micro-heterogeneous structure of the respective material. While the macroscopic deformation gradient $\overline{F}(\bar{x})$ defines the boundary conditions of an attached RVE as described in the latter section, the volumetric average over the microscopic Piola-Kirchhoff stresses in the RVE yield the macroscopic stress $\overline{P}(\overline{F}(\bar{x}))$. Therefore, in each macroscopic integration point $\bar{x}$, we obtain
\begin{equation}
	\overline{P}(\overline{F}(\bar{x})) = \frac{1}{V}\int_{B_0}P(F)dV
	\label{avg_stress}
\end{equation}
where $V=\lvert B_0 \rvert$ is the volume of the RVE belonging to $\bar{x}$. We now have to specify what the integral actually means in the case of an RVE modeled with beams. In that case, the integral over the RVE is computed as the sum over all beam elements. For each beam element, the integration of the first Piola-Kirchhoff stress is based on a procedure given in \cite{ARL_stress_av}; see below for further details. We thus obtain

\begin{equation*}
	\overline{P}(\overline{F}(\bar{x})) = \frac{1}{V}\sum_{e=1}^{n_B}\int_{V_e}P(F)dV.
\end{equation*}

Let us now describe the procedure of averaging the stresses within a single beam as given in \cite{ARL_stress_av}. We first consider a single beam element $e=[v_i,v_j]$ which is aligned with the unit vector $r=\left(r_x,r_y,r_z\right)^T=\frac{v_j-v_i}{L_e}$. Here, $L_e$ is the length of the considered beam element $e$. It is well-known that the first Piola-Kirchhoff stress can be expressed in terms of $P=J\,\sigma\,F^{-T},$ where $\sigma$ is the Cauchy stress tensor and $J$ is the determinant of the deformation gradient $F$. With the formulation of the locally transformed deformation $\tilde{u}$ as a cubic polynomial with respect to the normalized distance to the ends of the beam $\xi$, it is possible to derive the expressions $J=1+\text{tr}(\nabla u)$ and $F^{-1}=\frac{1}{J}\left((1+J)I-F\right)$. \Cref{appendix1} of the appendix presents a detailed derivation of these relations. Based on these equations we derive a formulation for the integral of the Piola-Kirchhoff stress over the volume of the beam element. This integral can be expressed as the sum of three separate integrals as

\begin{align}
	\begin{aligned}
		&\int_{V_e}PdV = \int_{V_e}\sigma \left(\left(1+J\right)I-F^{T}\right)dV	
		=&\int_{V_e}\sigma dV+\int_{V_e}J\sigma dV-\int_{V_e}\sigma F^{T} dV.
	\end{aligned}
	\label{eq:PK_integral}
\end{align}

The computation of the first integral of the right-hand side is carried out as in~\cite{ARL_stress_av}. There, a beam network is considered and a virtual work approach for the computation of an average stress for the network is introduced. The resulting expression for the average stress considers only the forces that appear in each beam element but not the higher moments. With this simplification the approach is not fully in alignment with the Euler-Bernoulli beam theory. However, it is to expect that the macroscopic behavior would not change significantly by the use of a more complex approach for the computation of the average stress. By applying this virtual work approach the first integral can be reformulated as
\begin{align}
	\begin{aligned}
		&\int_{V_e}\sigma dV = \int_{V_e}\left(\begin{smallmatrix}
			\sigma_{xx} & \sigma_{xy} & \sigma_{xz} \\
			\sigma_{xy} & \sigma_{yy} & \sigma_{yx} \\
			\sigma_{xz} & \sigma_{yz} & \sigma_{zz}
		\end{smallmatrix}\right) dV   \\
		=& \int_{V_e}L_e\left(\begin{smallmatrix}
			Q_x r_x & \frac{1}{2}(Q_x r_y+Q_y r_x) & \frac{1}{2}(Q_x r_z+Q_y r_z) \\
			\frac{1}{2}(Q_x r_y+Q_y r_x) & Q_y r_y & \frac{1}{2}(Q_y r_z+Q_z r_y)  \\
			\frac{1}{2}(Q_x r_z+Q_z r_x) & \frac{1}{2}(Q_y r_z+Q_z r_y) & Q_z r_z
		\end{smallmatrix}\right) dV
	\end{aligned}
\label{eq:sigma_integral}
\end{align}
with the forces $Q_x,Q_y,$ and $Q_z$ being expressed depending on the derivative of the projected deformation as $Q_{x_k}=-\frac{E\cdot I}{L_e^3}\cdot \frac{d^3 \tilde{u}_{x_k}}{d \xi^3}$.

Before deriving expressions for the other two integrals of the right-hand side of \Cref{eq:PK_integral}, we make the additional assumption that the deformation $u$ is constant over all cross sections of each beam element. This does also not align with the basic Euler-Bernoulli beam theory. However, this approximation is reasonable in our case since we are only considering relatively small cross sections compared to the length of each beam and additionally, we consider small deformations and therefore only a small bending occurs in the beams. The latter two assumptions do not need to be added, as they are already a requirement for the Euler-Bernoulli beam theory. Together with the equation derived for the integrated Cauchy stress we are able to formulate a simple expression for each of the remaining integrals. We finally obtain

\begin{align}
	\begin{aligned}
		&\int_{V_e} J\sigma \text{d}V = L_e\int_{0}^{1} J \text{d}\xi \, \int_{A_e} \sigma \, \text{d}A  \\
		=& \left(1+\frac{\langle u(v_2)-u(v_1), v_2-v_1\rangle}{L_e^2}\right)\int_{V_e} \sigma \, \text{d}V \; {\rm and}
	\end{aligned}
	\label{eq:Js_integral} \\
	\begin{aligned}
		&\int_{V_e}\sigma F^T\,\text{d}V = \int_{V_e}\sigma (I+(\nabla u)^T)\,\text{d}V  \\
		=& \int_{V_e}\sigma\,\text{d}V\cdot\left(I+\frac{\left(u(v_2)-u(v_1) \right)\cdot\left(v_2-v_1\right)^T}{L_e}\right).
	\end{aligned}
	\label{eq:sF_integral}
\end{align}

\begin{figure}
	\begin{minipage}{0.49\textwidth}
		\includegraphics[width=\textwidth]{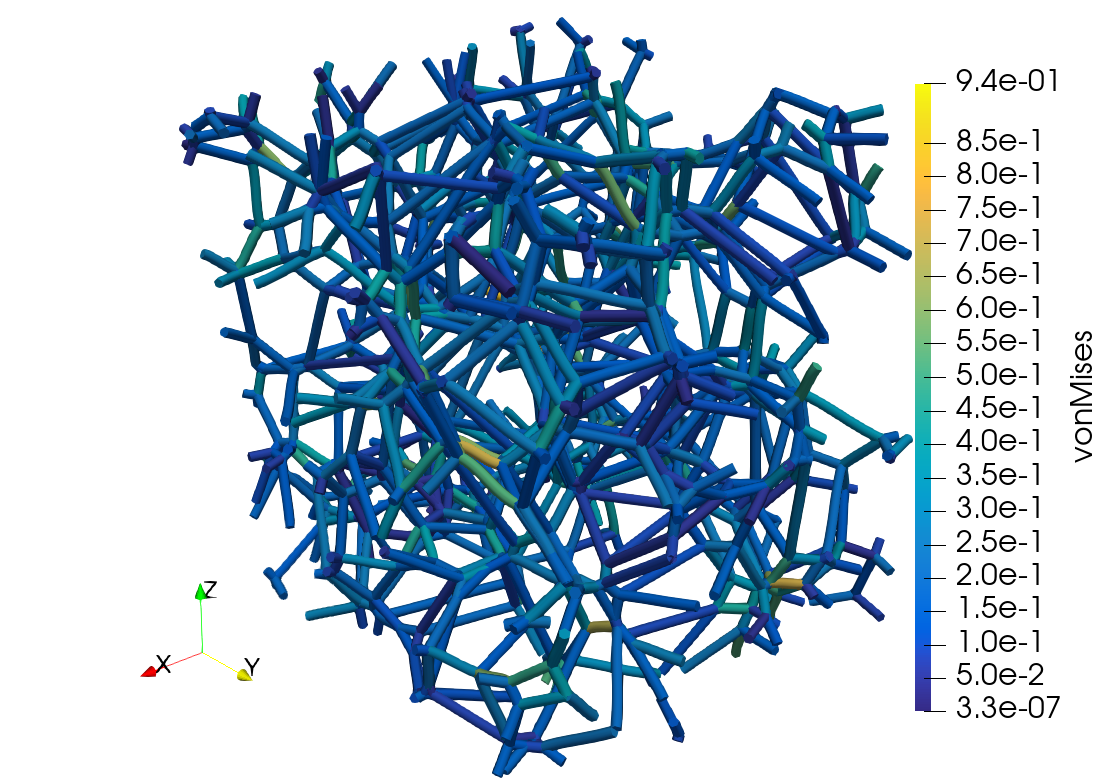}
	\end{minipage}
	\begin{minipage}{0.49\textwidth}
		$$ \overline{F} = \left(\begin{array}{ccc}
			0.9 & 0 & 0\\
			0 & 1 & 0 \\
			0 & 0 & 1
		\end{array}\right) $$
	{\huge$$\downarrow$$}
		$$\overline{P} = \left(\begin{array}{ccc}
			-111.53 & 0.53 & 0.43	\\
			0.53&-2.84& 1.79 	\\
			0.43&1.79& -7.25
		\end{array}\right)$$
	\end{minipage}
	\caption{Result of three-dimensional beam frame model that is compressed in $x$ direction. Image of the deformed RVE on the left and the deformation gradient with the resulting averaged first Piola-Kirchhoff stress tensor on the right. Colors of the beams represent the Von-Mises stress values in the corresponding elements. Von-Mises stress values are given in megapascal (MPa).}
	\label{fig:BF_micro_example}
\end{figure}

Due to the assumption that the deformation $u$ is constant over cross sections, it is possible to separate the volume integral into one integral over the surface of the cross section regarding the Cauchy stress $\sigma$ and one integral over the length of the beam regarding the  variables $J$ and $F$, which both depend on the deformation $u$. Summing \Cref{eq:sigma_integral}, \Cref{eq:Js_integral}, and \Cref{eq:sF_integral}, finally yields an expression for the complete integral on the left-hand side of \Cref{eq:PK_integral}. One example for the computation of the Piola-Kirchhoff stress tensor for one microscopic deformation is presented in \Cref{fig:BF_micro_example}.

\begin{figure*}
	\centering
	\includegraphics[width=0.7\textwidth]{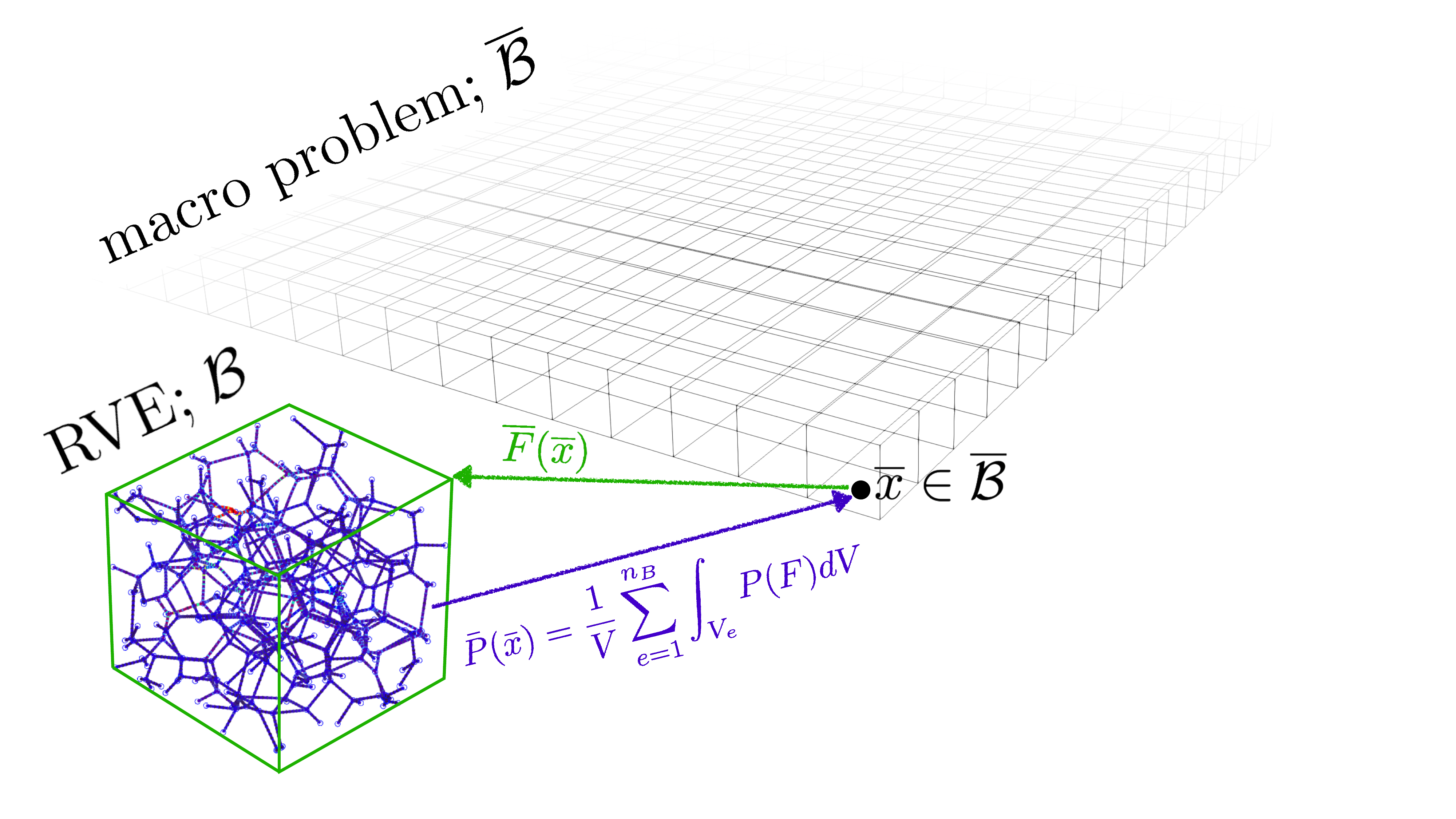}
	\caption{FE$^2$ method with an RVE modeled with beams. The macroscopic deformation gradient $\overline{F}(\bar{x})$ defines the boundary conditions of the beam frame model attached to the macroscopic integration point $\bar{x}$. The stresses in the beams after solving the RVE problem is averaged and results in the macroscopic stress $\overline{P}(\overline{F}(\bar{x}))$.}
	\label{fig:FE2_bf}
\end{figure*}

After computing the average stress in each beam, the average stress within the RVE $\overline{P}$ can be computed using \Cref{avg_stress} by summing over all beams. An algorithmic representation of our homogenization approach is shown in \Cref{fig:FE2_bf}.

\subsection{Algorithmic description}

The macroscopic problem given in \Cref{macro_problem}, which is discretized using finite elements, is solved using BFGS or an alternative (Quasi)-Newton method where the Jacobian is approximated using central difference quotients. For stability and robustness, a dynamic load stepping is integrated such that the total macroscopic deformation can be applied incrementally in several pseudo-time or load steps, respectively.  When integrating over the macroscopic finite elements using a Gauss quadrature rule, for each integration point the microscopic problem is solved to obtain the corresponding value of $\overline{P}$. Let us recapitulate that we use a linear beam frame model on the microscale and thus also the resulting macroscopic problem is actually linear. Using the exact Newton method as a macroscopic nonlinear solver would therefore lead to its convergence in a single step. However, since we have put no efforts in finding an exact formulation for a consistent tangent modulus $\frac{\partial \overline{P}(\overline{F}(\bar{x}))}{\partial \overline{F}(\bar{x})}$ we cannot use an exact Newton method. Instead, as already mentioned, we use the Quasi-Newton method BFGS and can avoid computing $\frac{\partial \overline{P}(\overline{F}(\bar{x}))}{\partial \overline{F}(\bar{x})}$ completely. As an alternative, we can approximate the Jacobian matrix using difference quotients, which is an alternative Quasi-Newton approach. For the approximation we use central differences of the form $$\left(\frac{\partial \overline{P}(\overline{F}(\bar{x}))}{\partial \overline{F}(\bar{x})}\right)_{i,j}\approx\frac{\overline{P}(\overline{F}(\bar{x})+\epsilon_j)_i-\overline{P}(\overline{F}(\bar{x})-\epsilon_j)_i}{2\cdot\epsilon}$$ for $i,j\in{1,\dots,4}$ or $i,j\in{1,\dots,9}$. For our applications the approximation has shown to yield robust results for $\epsilon=1e-6$. Let us remark that for the surrogate model introduced later on, we easily can use the exact Newton method and alternatively BFGS. Let us finally give a brief overview of the algorithm in \Cref{fig:algo}, where we describe the procedure within one load step of the dynamic load stepping procedure.

\begin{figure}
	\centering
	\begin{itemize}
	\item[] {\bf Init} macroscopic deformation $\bar{u}^{(0)}$ using current load
	\item[] {\it /*current load is defined by dynamic load stepping*/}
	\item[] {\bf Loop} over $k$ until convergence
	\begin{itemize}
	\item[] {\bf Loop} over all integration points $\bar{x}$
	\begin{itemize}
	\item[] {\bf Compute} $\overline{F}(\bar{x})$ from $\bar{u}^{(k)}$
	\item[] {\bf Apply} periodic boundary conditions to RVE defined by $\overline{F}(\bar{x})$
	\item[] {\bf Solve} RVE
	\item[] /*solving a beam frame problem*/
	\item[] {\bf Compute} $\overline{P}$ from RVE solution
	\end{itemize}
	\item[] {\bf EndLoop}
	\item[] {\bf Assemble} macroscopic residual vector using $\overline{P}$
	\item[] {\bf Compute} Quasi-Newton (e.g. BFGS) update $\delta \bar{u}^{(k)}$
	\item[] {\bf Update} $\bar{u}^{(k+1)} = \bar{u}^{(k)} + \delta \bar{u}^{(k)}$ 
	\item[] {\bf Check} for convergence of Quasi-Newton method 
	\end{itemize}
	\item[] {\bf EndLoop}
	\end{itemize}	
	\caption{FE$^2$ algorithm using BFGS or an alternative Quasi-Newton approach as a solver on the macroscopic scale and beam frame RVEs. The complete algorithm is usually embedded in a (dynamic) load stepping scheme for robustness and global convergence.}
	\label{fig:algo}
\end{figure}

\section{Neural Network-based Surrogate Model}\label{sec:NN}
Depending on the number of beam elements in the RVE, solving the microscopic problem and evaluating \Cref{avg_stress} can be computationally expensive. This becomes especially important for three-dimensional problems since a three-dimensional RVE naturally needs to consist of a much larger number of beam elements. Additionally six degrees of freedom instead of four are necessary, which also leads to larger systems of equations. To deal with this issue, we introduce a surrogate model that is supposed to compute the average Piola-Kirchhoff stress and replace the microscopic simulations based on the beam frame model. Here, we aim for training a machine learning-based surrogate model. More precisely, we have developed an artificial neural network (NN) which is trained to predict the average Piola-Kirchhoff stress tensor for one fixed microscopic structure. The NN has to be evaluated in each integration point of the macroscopic finite element problem and the input of the NN is always the deformation gradient $\overline{F}(\bar{x})$. The open-source package TensorFlow \cite{tensorflow2015-whitepaper} is used for the development and the training of each regarded NN in this section.

\subsection{Surrogate model in two dimensions}
Let us first describe the training procedure and give some details on the NN architecture we use. In a first step, the fixed microscopic beam frame problem is generated with the approach explained in \Cref{sec:methods}. We consider a two-dimensional open-porous structure and randomly set pores following a given pore size distribution. The resulting beam frame structure which is used to generate the training data for the model consists of 225 beam elements and 172 joints at which the beams are connected. Let us note that the model is comparably small and in the present article we just aim for a proof of concept of the suggested method. Also, we only use a prototype MATLAB implementation. It is planned to consider larger and more representative RVEs in the future using more efficient C/C++ based implementations.

For the two-dimensional RVE the data set for the training of the NN consists of about 34,000 pairs of macroscopic deformation gradients $\overline{F}$ and stresses $\overline{P}$, where the latter one is obtained by solving the corresponding RVE problem and averaging over the stresses within the beams. To obtain a large variety of different macroscopic deformation gradients, several FE$^2$ simulations following the algorithm from \Cref{fig:algo} have been carried out and for each BFGS step in all integration points the pairs $(\overline{F},\overline{P})$ are stored.

\begin{figure*}
	\begin{minipage}{0.24\textwidth}
		\includegraphics[width=\textwidth]{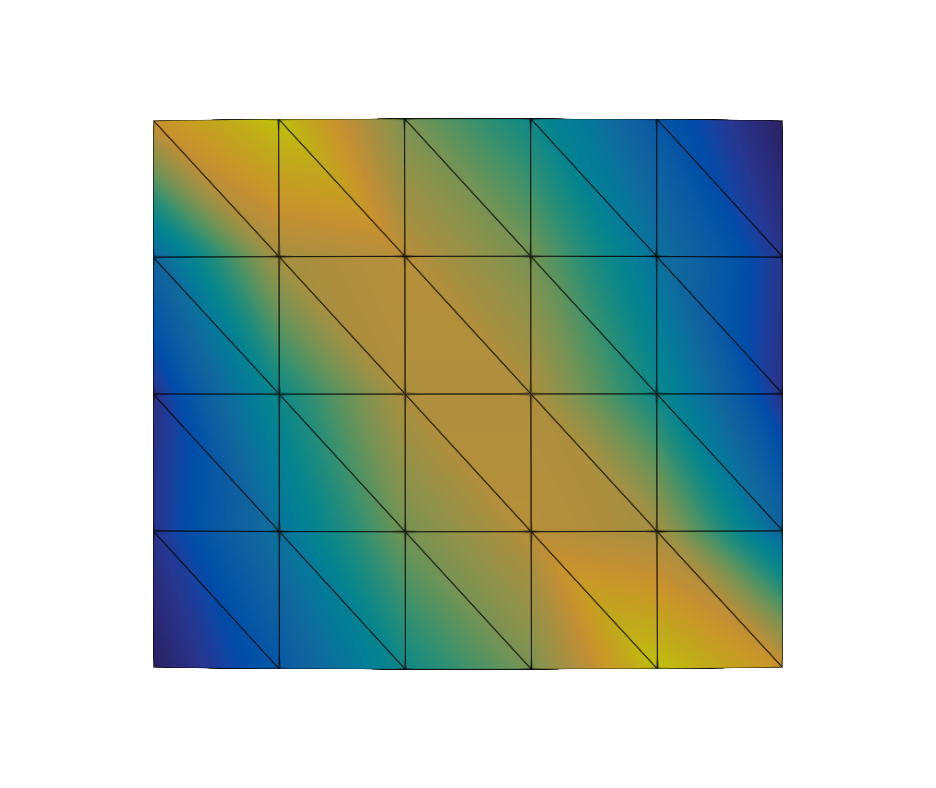}
	\end{minipage}
	\begin{minipage}{0.24\textwidth}
		\includegraphics[width=\textwidth]{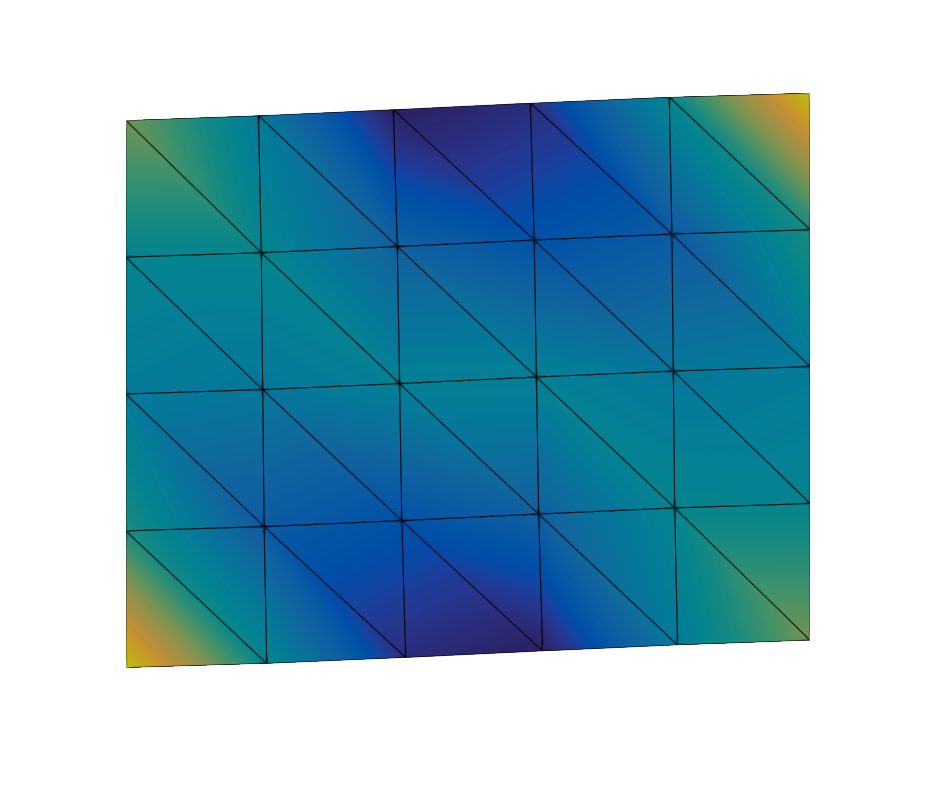}
	\end{minipage}
	\begin{minipage}{0.24\textwidth}
		\includegraphics[width=\textwidth]{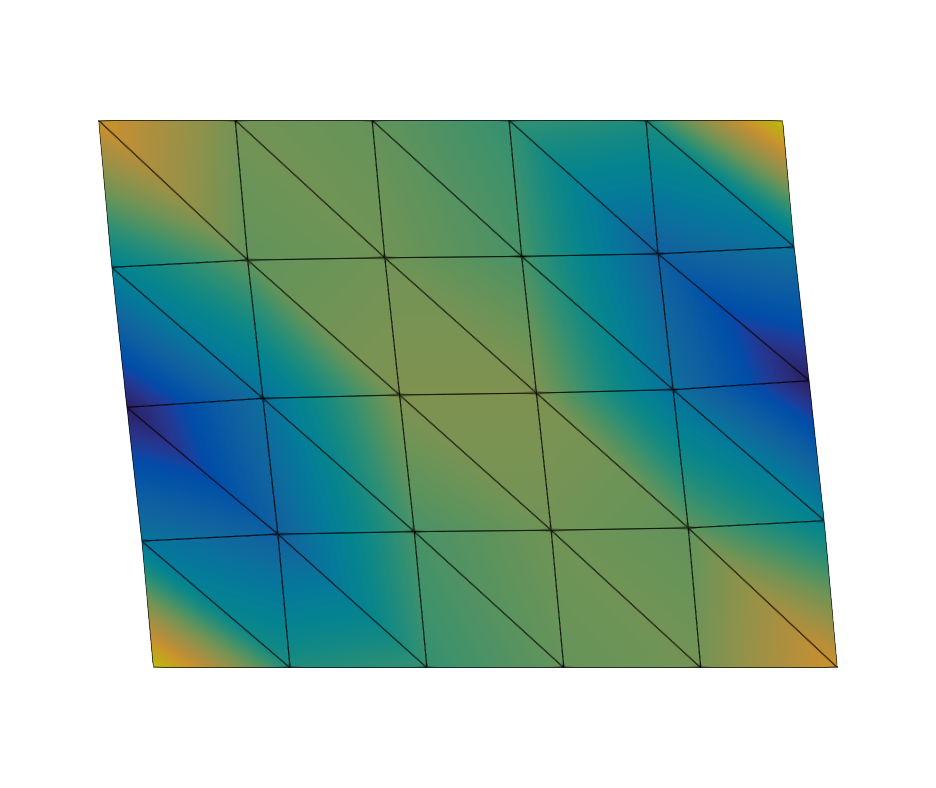}
	\end{minipage}
	\begin{minipage}{0.24\textwidth}
		\includegraphics[width=\textwidth]{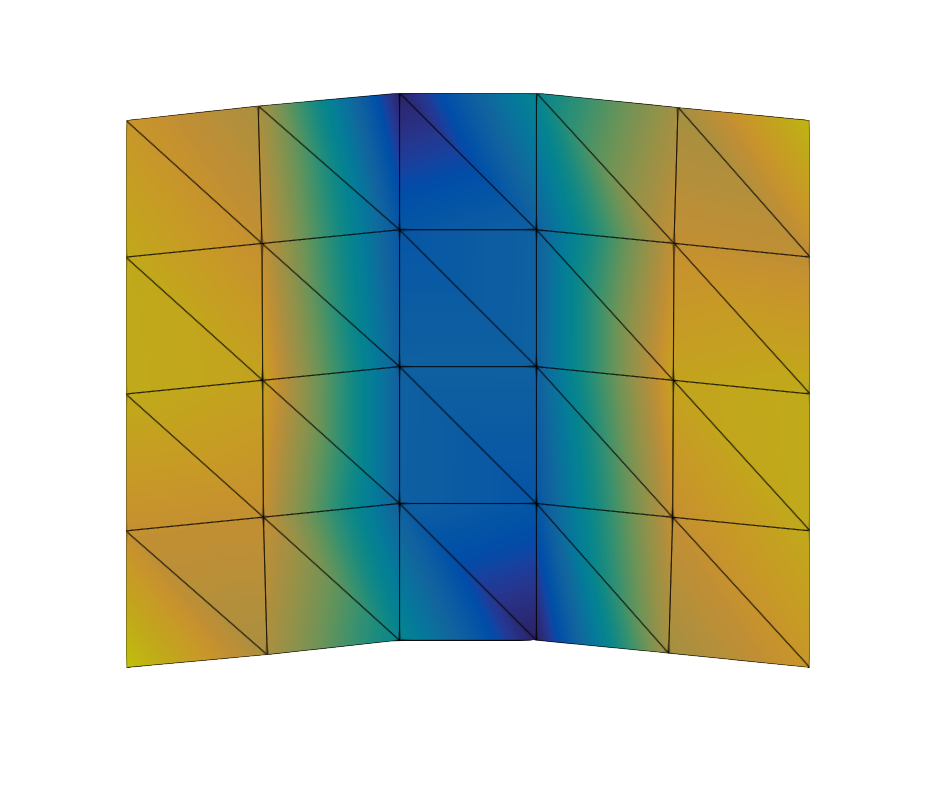}
	\end{minipage}
	
	\begin{minipage}{0.12\textwidth}
		\,
	\end{minipage}
	\begin{minipage}{0.24\textwidth}
		\includegraphics[width=\textwidth]{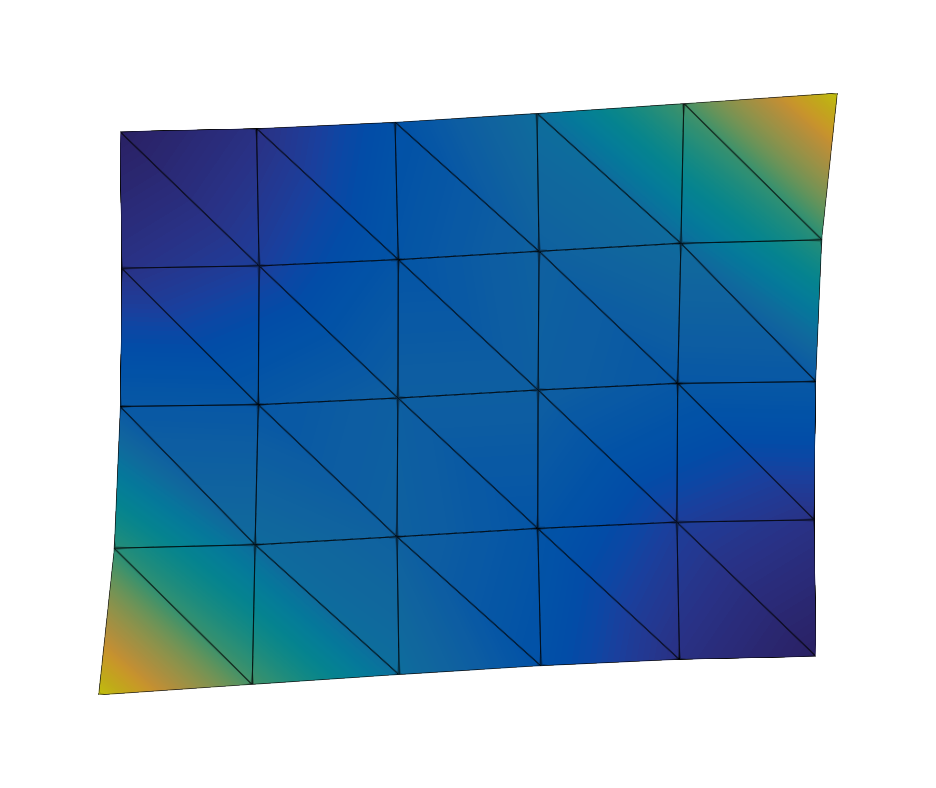}
	\end{minipage}
	\begin{minipage}{0.24\textwidth}
		\includegraphics[width=\textwidth]{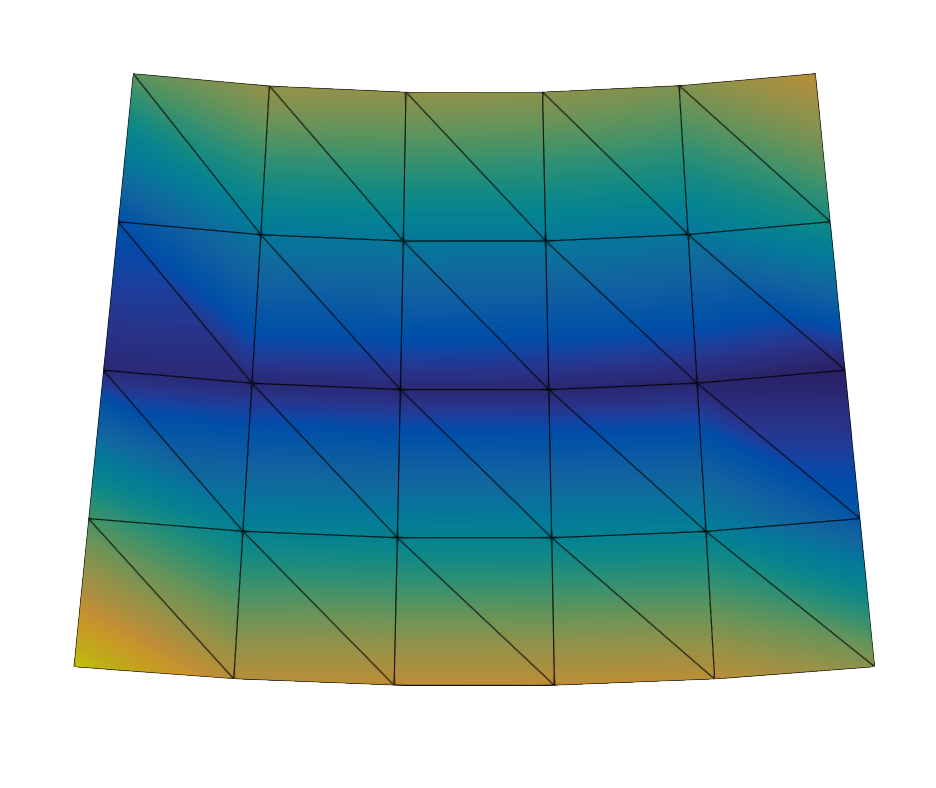}
	\end{minipage}
	\begin{minipage}{0.24\textwidth}
		\includegraphics[width=\textwidth]{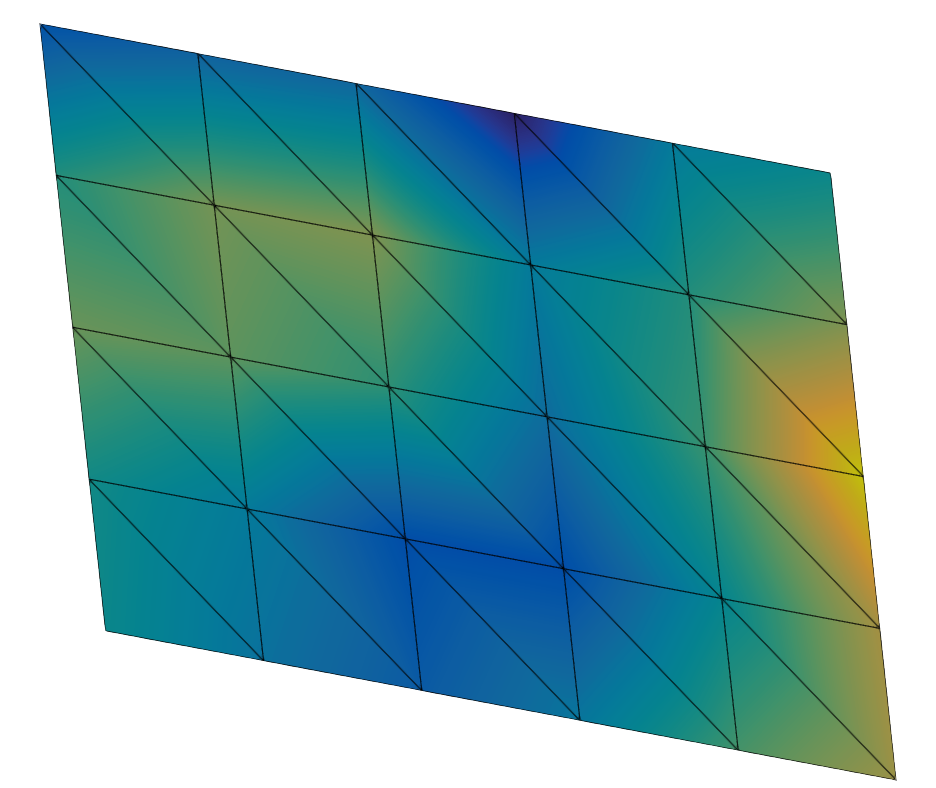}
	\end{minipage}
	\caption{Seven basic deformations that are used to generate the training data in two dimensions. To obtain a large training data set, the basic deformations are modified randomly by changing the grade and orientation of deformation.}
	\label{fig:train_examples2d}
\end{figure*}

For setting up the different FE$^2$ simulations we defined seven different macroscopic basic test cases and further vary these by changing the degree of deformation for each case. Each of the seven basic test cases has different macroscopic boundary conditions and a different deformation of the material such that the data set is expected to cover a wide variety of microscopic deformations. The magnitude and direction of each computed macroscopic deformation is set randomly during the generation process. Results for the deformation examples that are used for the generation are presented in \Cref{fig:train_examples2d}. The same relatively small grid is used for each of the examples and the resulting problem has 60 degrees of freedom. There is no need to use larger problems for the generation of the training data. In contrast, the computation of a larger number of small problems yields more diversity among the localized deformation gradients in the integration points compared with a smaller number of large problems. Therefore, we expect a higher variety in the generated training data defining and using many small macroscopic problems.

In the generation of the training data we have exclusively used linear finite elements to discretize the macroscopic problem. However, the choices of the basis functions and elements for the macroscopic finite element method are not expected to affect the quality of the training data since the generated input and output values of the NN are only related to the microscopic problem. The validation data set consisting of about 4,000 input and output pairs has been generated in the same manner.

\begin{figure}
	\begin{minipage}{0.49\textwidth}
		\includegraphics[width=\textwidth]{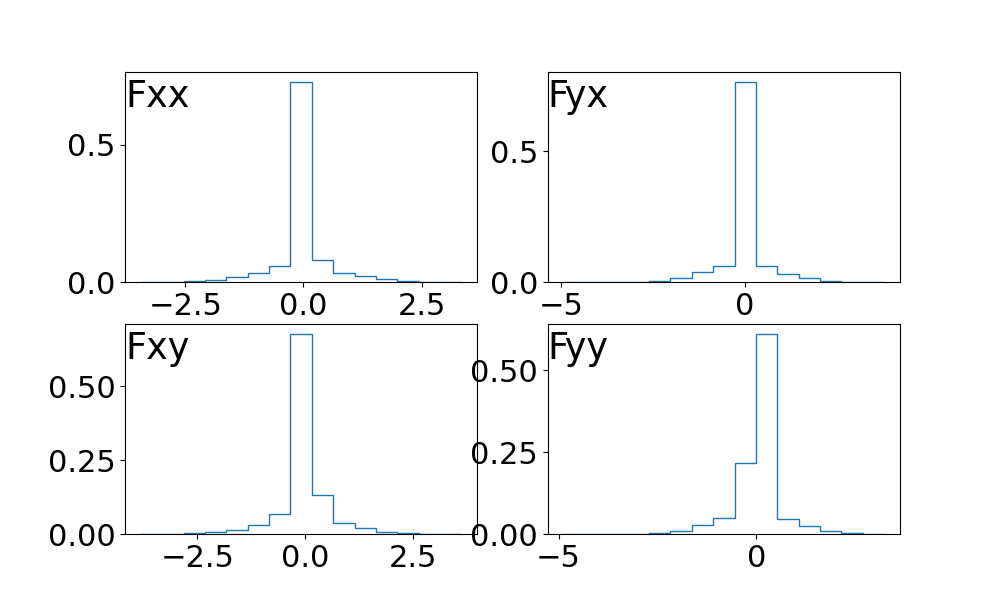}
	\end{minipage}
	\begin{minipage}{0.49\textwidth}
		\includegraphics[width=\textwidth]{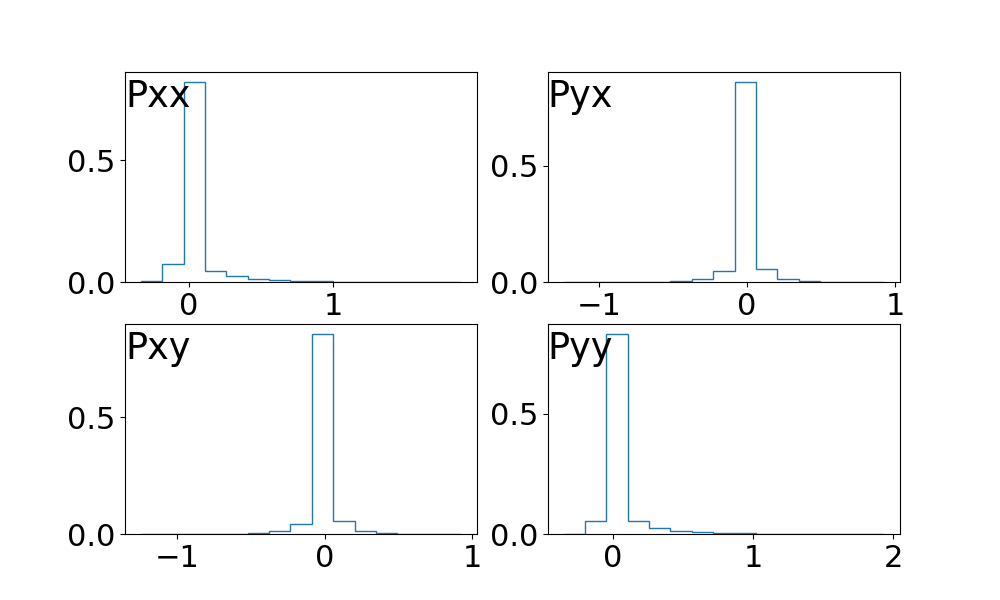}
	\end{minipage}
	\caption{Distribution of the input and output data generated for the training of the two-dimensional neural network.}
	\label{fig:data_disrtib2D}
\end{figure}

The distribution of the generated data set is presented in \Cref{fig:data_disrtib2D}. All input components are distributed relatively tight around zero. The corresponding small values for the standard deviations which can be observed in \Cref{fig:input_data_tab2D} are intended since we are only assuming small deformations. This assumption is necessary because the applied beam frame model which is based on Euler-Bernoulli beam theory is only useful for small deformations. This is acceptable because a large proportion of reported literature investigates the linear elastic properties of open-porous materials owing to their dependence on the density that is demonstrated in terms of scaling laws, particularly between Young's modulus or compressive strength versus density \cite{gibson-ashby-1999, aney-jsst-2023}.

As illustrated in \Cref{fig:output_data_tab2D}, the components of the resulting first Piola-Kirchhoff stress tensors which are the target output variables for the training of our neural network are similarly distributed around the mean of zero. This is not surprising since the matrix in \Cref{fig:correlaion2D} presents distinct correlations between the input and output components.

\begin{table}
	\begin{tabular}{|c|cccc|}
		\hline
		variable name & $\overline{F}_{xx}$ & $\overline{F}_{yx}$ & $\overline{F}_{xy}$ & $\overline{F}_{yy}$	\\
		\hline
		mean & -0.01 & 0 & -0.01 & -0.01 \\
		standard deviation & 0.52 & 0.55 & 0.53 & 0.54 \\
		minimum value & -3.44 & -4.93 & -3.81 & -4.87	\\
		maximum value & 3.34 & 3.80 & 3.63 & 3.23	\\
		\hline
	\end{tabular}
	\caption{Description of the input data for the training of the two-dimensional neural network}
	\label{fig:input_data_tab2D}
\end{table}

\begin{table}
	\begin{tabular}{|c|cccc|}
		\hline
		variable name & $\overline{P}_{xx}$ & $\overline{P}_{yx}$ & $\overline{P}_{xy}$ & $\overline{P}_{yy}$	\\
		\hline
		mean & 0.03 & 0 & 0 & 0.04 \\
		standard deviation & 0.14 & 0.09 & 0.09 & 0.15 \\
		minimum value & -0.33 & -1.23 & -1.25 & -0.35	\\
		maximum value & 1.88 & 0.93 & 0.93 & 1.93	\\
		\hline
	\end{tabular}
	\caption{Description of the output data for the training of the two-dimensional neural network.}
	\label{fig:output_data_tab2D}
\end{table}

\begin{table}
	\begin{tabular}{|c|cccc|}
		\hline
		& $\overline{P}_{xx}$ & $\overline{P}_{yx}$ & $\overline{P}_{xy}$ & $\overline{P}_{yy}$	\\
		\hline
		$\overline{F}_{xx}$ & \colorbox{orange}{0.66} & -0.01 & -0.01 & 0.04 \\
		$\overline{F}_{yx}$ & 0 & \colorbox{yellow}{0.50} & \colorbox{yellow}{0.50} & -0.02 \\
		$\overline{F}_{xy}$ & 0.01 & \colorbox{yellow}{0.52} & \colorbox{yellow}{0.52} & 0.01	\\
		$\overline{F}_{yy}$ & 0.05 & -0.01 & -0.01 & \colorbox{orange}{0.62}	\\
		\hline
	\end{tabular}
	\caption{Correlation coefficients between the input and output variables for the two-dimensional data set.}
	\label{fig:correlaion2D}
\end{table}

\begin{figure}[h!]
\begin{center}
	\begin{tikzpicture}[scale=0.53]
		\def\height{0}
		\def\depth{0}
		\def\width{2.5}
		
		\def\nfirstlayer{4}
		\def\nfifthlayer{4}

		\node at (-0.8,\height,\depth) {$\mathbf{\overline{F}}$};
		\node at (4*\width+0.8,\height,\depth) {$\mathbf{\overline{P}}$};
		
		\foreach \x in {1,...,\nfirstlayer}
		{
			\draw[thick, blue] (0, \height+\nfirstlayer/2+0.5-\x, \depth) circle (6pt);
			\foreach \y in {1,...,2}
			{
				\draw[-stealth,line width=1pt,gray] (0, \height+\nfirstlayer/2+0.5-\x, \depth) -- (\width, \height+\y, \depth);
				\draw[-stealth,line width=1pt,gray] (0, \height+\nfirstlayer/2+0.5-\x, \depth) -- (\width, \height-\y, \depth);
			}
		}
		\foreach \x in {1,...,2}
		{
			\draw[thick, blue] (\width, \height+\x, \depth) circle (6pt);
			\draw[thick, blue] (\width, \height-\x, \depth) circle (6pt);
			\foreach \y in {1,...,2}
			{
				\draw[-stealth,line width=1pt,gray] (\width, \height+\x, \depth) -- (2*\width, \height+\y, \depth);
				\draw[-stealth,line width=1pt,gray] (\width, \height+\x, \depth) -- (2*\width, \height-\y, \depth);
				\draw[-stealth,line width=1pt,gray] (\width, \height-\x, \depth) -- (2*\width, \height+\y, \depth);
				\draw[-stealth,line width=1pt,gray] (\width, \height-\x, \depth) -- (2*\width, \height-\y, \depth);
			}
		}
		\node at (\width,\height,\depth) {\textbf{\vdots}};
		\foreach \x in {1,...,2}
		{
			\draw[thick, blue] (2*\width, \height+\x, \depth) circle (6pt);
			\draw[thick, blue] (2*\width, \height-\x, \depth) circle (6pt);
			\foreach \y in {1,...,2}
			{
				\draw[-stealth,line width=1pt,gray] (2*\width, \height+\x, \depth) -- (3*\width, \height+\y, \depth);
				\draw[-stealth,line width=1pt,gray] (2*\width, \height+\x, \depth) -- (3*\width, \height-\y, \depth);
				\draw[-stealth,line width=1pt,gray] (2*\width, \height-\x, \depth) -- (3*\width, \height+\y, \depth);
				\draw[-stealth,line width=1pt,gray] (2*\width, \height-\x, \depth) -- (3*\width, \height-\y, \depth);
			}
		}
		\node at (2*\width,\height,\depth) {\textbf{\vdots}};
		\foreach \x in {1,...,2}
		{
			\draw[thick, blue] (3*\width, \height+\x, \depth) circle (6pt);
			\draw[thick, blue] (3*\width, \height-\x, \depth) circle (6pt);
			\foreach \y in {1,...,\nfifthlayer}
			{
				\draw[-stealth,line width=1pt,gray] (3*\width, \height+\x, \depth) -- (4*\width, \height+\nfifthlayer/2+0.5-\y, \depth);
				\draw[-stealth,line width=1pt,gray] (3*\width, \height-\x, \depth) -- (4*\width, \height+\nfifthlayer/2+0.5-\y, \depth);
			}
		}
		\node at (3*\width,\height,\depth) {\textbf{\vdots}};
		\foreach \x in {1,...,\nfifthlayer}
		{
			\draw[thick, blue] (4*\width, \height+\nfifthlayer/2+0.5-\x, \depth) circle (6pt);
		}
		\node at (0,\height-2.7,\depth) {\#neurons:};
		\node at (\width,\height-2.7,\depth) {\textbf{128}};
		\node at (2*\width,\height-2.7,\depth) {\textbf{256}};
		\node at (3*\width,\height-2.7,\depth) {\textbf{128}};
		
		\node at (0.5*\width,\height+2.4,\depth) {\includegraphics[width=5mm]{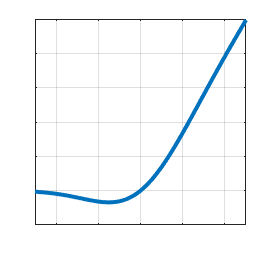}};
		\node at (1.5*\width,\height+2.6,\depth) {\includegraphics[width=5mm]{images/gelu.png}};
		\node at (2.5*\width,\height+2.6,\depth) {\includegraphics[width=5mm]{images/gelu.png}};
		\node at (3.5*\width,\height+2.4,\depth) {\includegraphics[width=5mm]{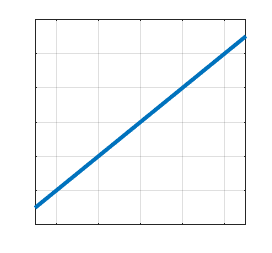}};
	\end{tikzpicture}
\end{center}
\caption{Architecture of the feed-forward neural network for a two-dimensional microscopic problem. First three layers have GELU activation and linear activation function is used for the output layer.}
\label{fig:NN}
\end{figure}

For the two-dimensional RVE a feed-forward NN with three hidden layers is used as the basic structure of our model (see \Cref{fig:NN}). The activation function for each hidden layer is chosen to be the Gaussian error linear unit (GELU) activation function \cite{hendrycks2023gaussian} and the number of neurons in the hidden layers are set to 128, 256, and 128. The activation of the output is linear. For choosing a proper architecture for the NN multiple activation functions in combination with different layer sizes were considered. A grid search evaluation covering these parameters has shown that GELU activation suits especially well for the given problem and that a higher number of neurons in the model can lead to a smaller training error. However, the reduction of the error is only significant up to 256 neurons per layer. The grid search results are presented in \Cref{fig:grid_search2d_tab} and the smallest values in the resulting training loss are highlighted. The table shows that a network architecture of 3 layers with 128 or 256 neurons per layer yields the lowest loss. Further testing with these parameters has shown that the selected architecture with 128, 256, 128 neurons is able to slightly reduce the loss even further.

\begin{table}
\begin{tabular}{|lc|ccc|}
	\hline
	activation& neur./layer& one layer & two layers & three layers	\\
	\hline
	sigmoid&64 &6.51e-05 & 1.76e-05 & 8.37e-06 \\
	&128&6.22e-05 & 7.64e-06 & 2.50e-06 \\
	&256&6.49e-05 & 3.91e-06 & 7.54e-07 \\
	&512&6.92e-05 & 8.54e-07 & 5.39e-07 \\	\hline
	tanh&64&2.98e-05 & 3.99e-06 & 2.13e-06 \\
	&128&2.62e-05 & 2.36e-06 & 4.06e-07 \\
	&256&2.72e-05 & 1.52e-06 & 2.32e-07 \\
	&512&2.83e-05 & 1.20e-06 & 2.38e-07 \\	\hline
	gelu&64&7.26e-07 & 1.27e-07 & 2.44e-08 \\
	&128&3.46e-07 & 3.65e-08 & \colorbox{yellow}{7.97e-09} \\
	&256&1.28e-07 & 1.72e-08 & \colorbox{yellow}{6.41e-09} \\
	&512&2.57e-07 & 1.11e-08 & 1.14e-08 \\
	\hline
\end{tabular}
\caption{Grid search results (training loss) for the neural network in two dimensions with different numbers of layers and neurons per layer.}
\label{fig:grid_search2d_tab}
\end{table}

\begin{figure}[]
	\centering
	\includegraphics[width=0.5\textwidth]{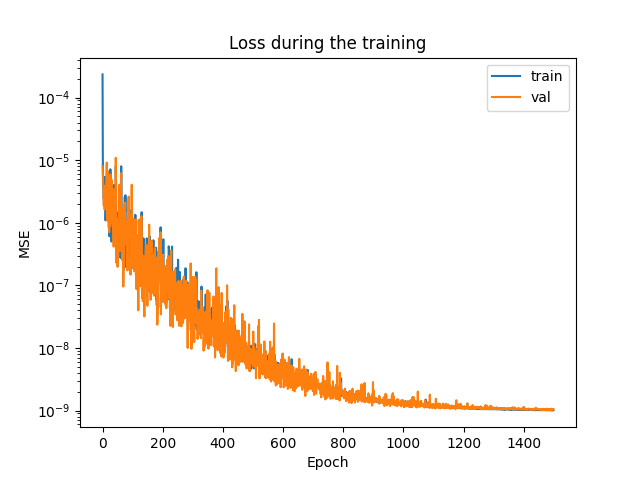}
	\caption{Loss during the training of the neural network for two-dimensional RVEs.}
	\label{fig:NNloss}
\end{figure}

For training, an adam optimizer \cite{kingma2017adam} has been used to minimize the mean squared error (MSE). After a total of 1,500 training epochs a sufficient reduction of the loss has been reached; see \Cref{fig:NNloss}. The loss for the validation data has also decreased sufficiently.

\subsection{Surrogate model in three dimensions}
For a given three-dimensional aerogel RVE the general procedure of generating training data and training the NN works similar to the two-dimensional case. The size of the input and output of the NN is different since the three-dimensional deformation gradient and Piola-Kirchhoff stress tensor have nine components each. We consider a microscopic RVE of 1,482 beam elements which are connected at 1,005 joints.

The training data is similarly generated from solving macroscopic deformation test cases. However, for the three-dimensional model the approach on how to set up the macroscopic tests differs since we do not use a fixed number of deformation examples. For the generation of the data set a relatively small cube geometry is considered with 192 degrees of freedom. Dirichlet boundary conditions are applied to each of the nodes on the boundaries of this geometry with a fixed deformation determined by a randomly generated deformation gradient. This means that the deformation on the boundaries is determined by $\bar{u}=\overline{F}_d \cdot \bar{x}$ with the matrix $\overline{F}_d$ being randomly generated. Due to the large influence of randomness, the procedure is expected to yield a high variety of deformations. For the training data about 114,000 data points are generated. The validation data consists of about 13,000 equally generated data points.

The input variables distributed around zero similar to the two-dimensional data set. The differed procedure for generating the training data results in a smaller difference between highest and lowest value and also small standard deviations of the variables. The data distribution can be visually observed in \Cref{fig:data_disrtib3D} and the values are also presented in \Cref{appendix2} of the appendix.

\begin{figure}
	\begin{minipage}{0.49\textwidth}
		\includegraphics[width=\textwidth]{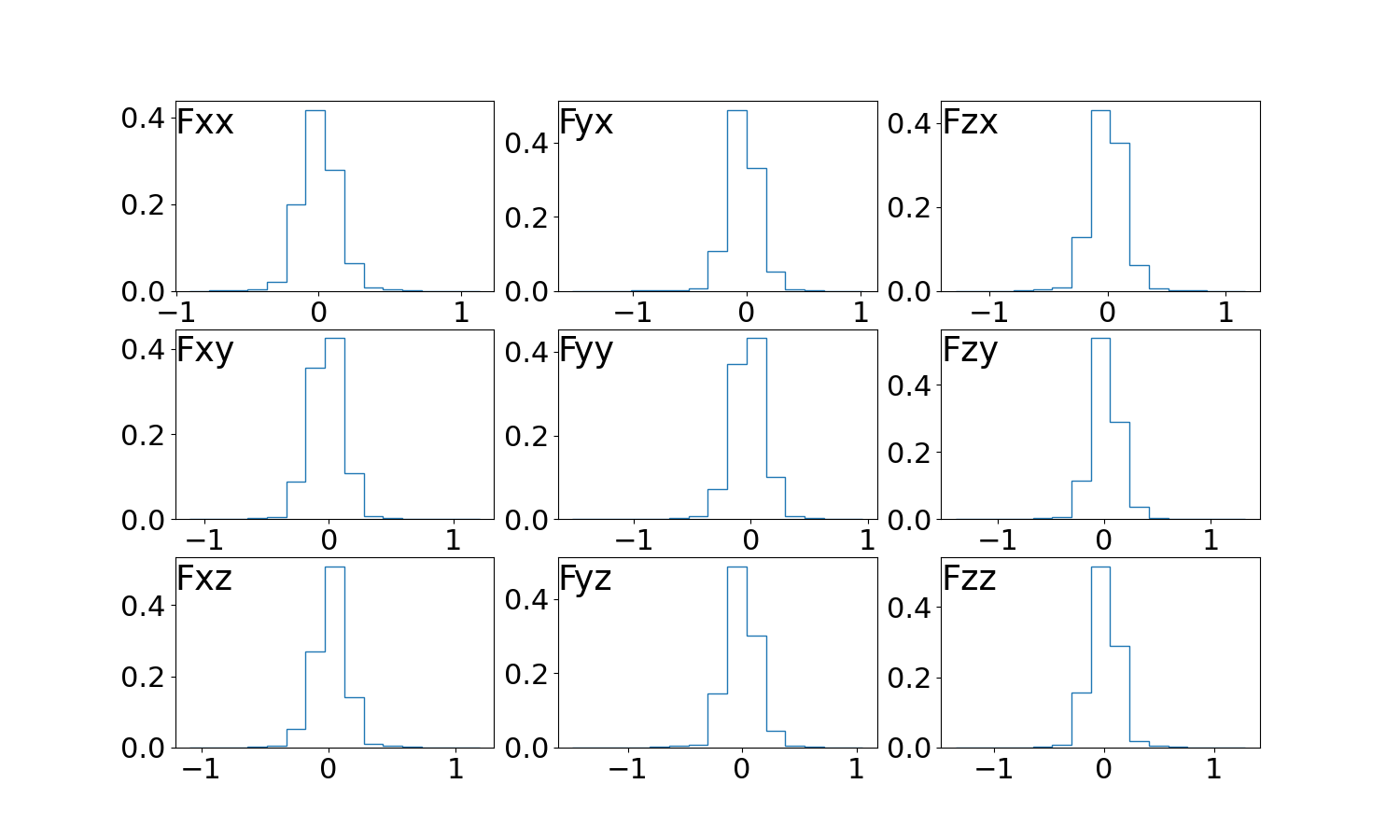}
	\end{minipage}
	\begin{minipage}{0.49\textwidth}
		\includegraphics[width=\textwidth]{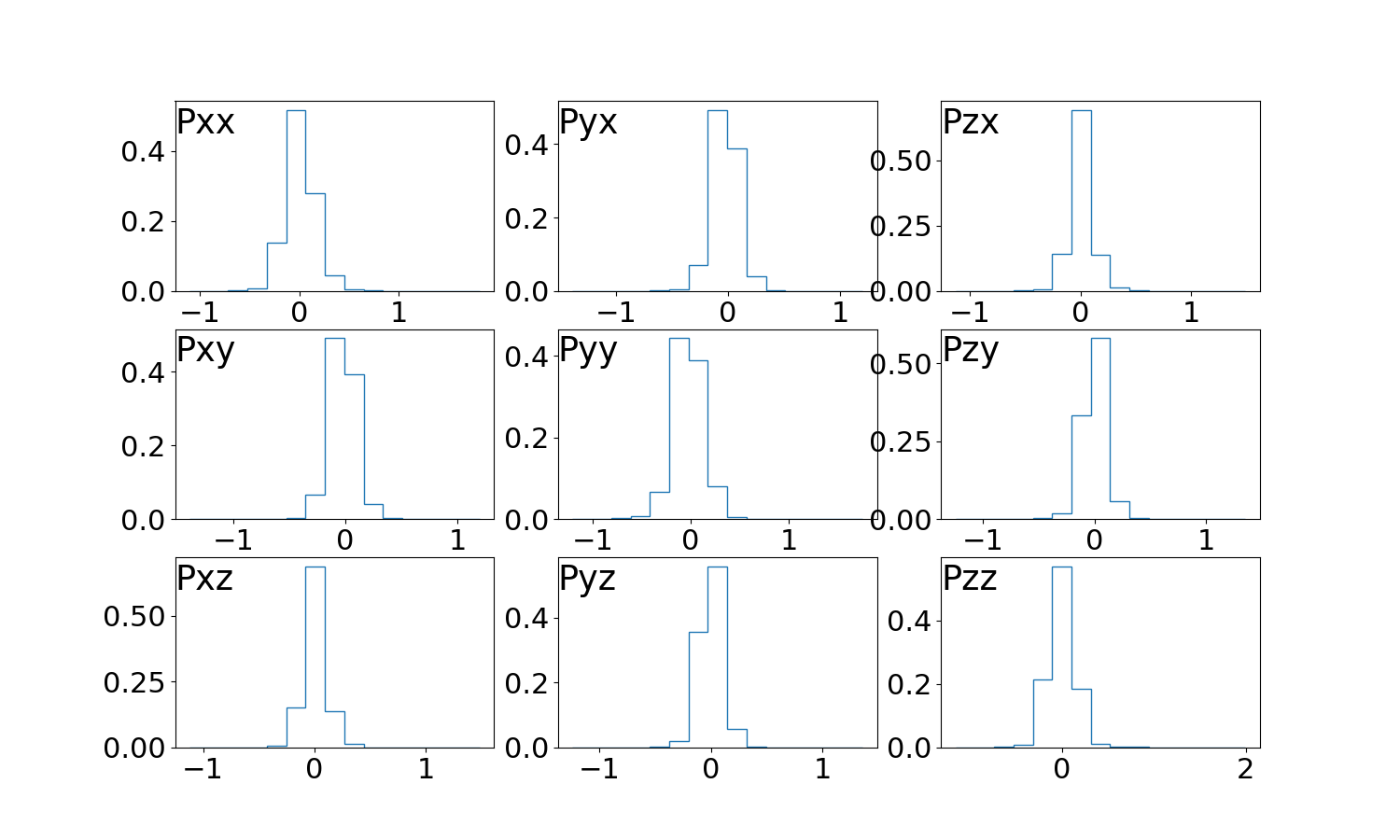}
	\end{minipage}
	\caption{Distribution of the input and output data generated for the training of the three-dimensional neural network.}
	\label{fig:data_disrtib3D}
\end{figure}

\begin{table}
	\begin{tabular}{|c|ccccccccc|}
		\hline
		& $\overline{P}_{xx}$ & $\overline{P}_{yx}$ & $\overline{P}_{zx}$ & $\overline{P}_{xy}$ & $\overline{P}_{yy}$ & $\overline{P}_{zy}$ & $\overline{P}_{xz}$ & $\overline{P}_{yz}$ & $\overline{P}_{zz}$ \\
		\hline
		$\overline{F}_{xx}$ & \colorbox{orange}{0.99} & -0.06 & -0.13 & -0.06 & 0.09 & -0.07 & -0.13 & -0.07 & 0.04 \\
		$\overline{F}_{yx}$ & 0.05 & \colorbox{yellow}{0.74} & 0.08 & \colorbox{yellow}{0.73} & -0.06 & -0.05 & 0.08 & -0.05 & -0.04 \\
		$\overline{F}_{zx}$ & -0.10 & 0.05 & \colorbox{yellow}{0.70} & 0.05 & 0.04 & 0.07 & \colorbox{yellow}{0.69} & 0.07 & 0.02 \\
		$\overline{F}_{xy}$ & -0.11 & \colorbox{yellow}{0.72} & -0.12 & \colorbox{yellow}{0.72} & 0.00 & -0.07 & -0.12 & -0.07 & -0.09 \\
		$\overline{F}_{yy}$ & 0.12 & -0.03 & 0.15 & -0.03 & \colorbox{orange}{0.99} & 0.02 & 0.15 & 0.02 & 0.07 \\
		$\overline{F}_{zy}$ & -0.01 & -0.11 & 0.00 & -0.11 & 0.08 & \colorbox{yellow}{0.71} & 0.00 & \colorbox{yellow}{0.71} & -0.01 \\
		$\overline{F}_{xz}$ & -0.08 & -0.10 & \colorbox{yellow}{0.66} & -0.10 & 0.16 & 0.02 & \colorbox{yellow}{0.66} & 0.02 & 0.06 \\
		$\overline{F}_{yz}$ & -0.05 & 0.03 & 0.09 & 0.03 & -0.07 & \colorbox{yellow}{0.70} & 0.09 & \colorbox{yellow}{0.70} & 0.13  \\
		$\overline{F}_{zz}$ & 0.04 & -0.09 & 0.06 & -0.09 & 0.03 & 0.06 & 0.06 & 0.06 & \colorbox{orange}{0.98} \\
		\hline
	\end{tabular}
	\caption{Correlation coefficients between the input and output variables for the three-dimensional data set}
	\label{fig:correlaion3D}
\end{table}

Besides the number of neurons in the input and output layer, the network architecture for the three-dimensional case is quite similar to the two-dimensional case. Here, the NN has two instead of three hidden layers with 512 neurons in each layer and again uses a GELU activation for each of these hidden layers. The choice of the architecture resulted again from a grid search evaluation. In \Cref{fig:grid_search3d_tab} the resulting training losses of the grid search algorithm is presented with the smallest values being highlighted. The grid search evaluation shows similar results to the two-dimensional case with the GELU activation function yielding the smallest loss values.

\begin{table}
	\begin{tabular}{|lc|ccc|}
		\hline
		activation& neur./layer& one layer & two layers & three layers	\\
		\hline
		sigmoid&64 &9.55e-06 & 3.75e-06 & 4.08e-06 \\
		&128&1.02e-05 & 1.96e-06 & 1.50e-06 \\
		&256&9.83e-06 & 7.15e-07 & 7.10e-07 \\
		&512& 1.17e-05 & 3.75e-07 & 7.72e-07 \\
		\hline
		tanh&64&3.31e-06 & 7.21e-07 & 3.80e-07 \\
		&128&2.94e-06 & 3.92e-07 & 1.81e-07 \\
		&256 &2.88e-06 & 4.01e-07 & 1.26e-07 \\
		&512 &3.03e-06 & 2.78e-07 & 1.62e-07 \\
		\hline
		gelu&64 &1.94e-07 & 3.53e-08 & 3.22e-08 \\
		&128&1.11e-07 & 1.21e-08 & 1.26e-08 \\
		&256&5.26e-08 & \colorbox{yellow}{5.46e-09} & 6.46e-09 \\
		&512 &3.70e-08 & \colorbox{yellow}{5.16e-09} & \colorbox{yellow}{5.15e-09} \\
		\hline
	\end{tabular}
	\caption{Grid search results (training loss) for the neural network in three dimensions with different numbers of layers and neurons per layer.}
	\label{fig:grid_search3d_tab}
\end{table}

\begin{figure}
	\begin{center}
		\includegraphics[width=0.5\textwidth]{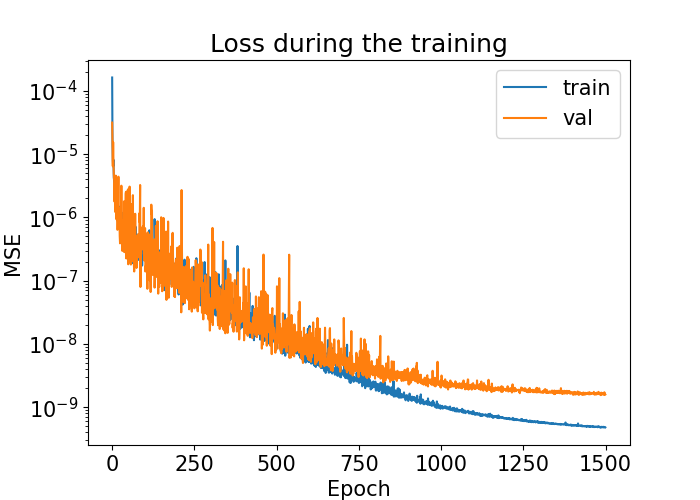}
	\end{center}
	\caption{Loss during the training of neural network for three-dimensional RVEs.}
	\label{fig:NNloss3d}
\end{figure}

During the training of the selected neural network architecture the training loss reduced within 1500 epochs to a final loss of 4.8e-10. The development of the loss during training is presented in \Cref{fig:NNloss3d}. As the image also shows, the validation loss has also decreased during the training though it settles on a slightly higher level compared to the training loss.

\subsection{Reduction of computational effort}
The surrogate model has been originally built to reduce the computational effort for computing the first Piola-Kirchhoff stress. Solving the beam frame model and computing the homogenized Piola-Kirchhoff stress tensor for the two-dimensional RVE takes about 200 times longer than the evaluation of the NN. It needs to be mentioned that these computation times have been measured with \textbf{MATLAB} and that a more efficient implementation with a better use of the hardware could yield another ratio. However, if we only consider the time  to solve the beam frame model's linear system of equations with \textbf{MATLAB}'s backslash operator, we still reach a reduction factor of around 100 with the surrogate model. Since the backslash operator, that is, using the efficient UMFPACK factorization package, is generally regarded as a fast option for solving a linear system of equations, it is expected that with another implementation the computational effort of the NN solver would still be at least 100 times smaller than the effort of the beam frame solver for this specific example.

For the three-dimensional case the factor is much higher since the beam frame RVE consists of more beam elements and solving the resulting system requires a higher computational effort. In this case, the evaluation of the NN takes about 7000 times less time than solving of the beam frame model. Here, the time needed in MATLAB for averaging the resulting stresses is included. Measuring only the time for applying UMFPACK, the NN is still about 4000 times faster. Let us remark that in realistic applications larger RVEs will be needed. Since the computing time of sparse direct solvers as UMFPACK grows approximately cubically with the dimension of the problem, we expect an even larger benefit using NN-based surrogate models.

\subsection{Algorithmic description}
Basically, the resulting algorithm exploiting the surrogate model is similar to algoritm \ref{fig:algo}. Only the evaluation of the average stress $\overline{P}$ is replaced by an evaluation of the NN. Let us remark that by deriving the NN, we can easily compute the tangent/Jacobian matrix and alternatively use Newton's method instead of BFGS here. We summarize the algorithm in \Cref{fig:algo2}.
\begin{figure}
	\centering
	\begin{itemize}
	\item[] {\bf Init} macroscopic deformation $\bar{u}^{(0)}$ using current load
	\item[] {\it /*current load is defined by dynamic load stepping*/}
	\item[] {\bf Loop} over $k$ until convergence
	\begin{itemize}
	\item[] {\bf Loop} over all integration points $\bar{x}$
	\begin{itemize}
	\item[] {\bf Compute} $\overline{F}(\bar{x})$ using $\bar{u}^{(k)}$
	\item[] {\bf Evaluate} NN with input $\overline{F}(\bar{x})$ to obtain output $\overline{P}$
	\end{itemize}
	\item[] {\bf EndLoop}
	\item[] {\bf Assemble} macroscopic residual vector using $\overline{P}(\overline{F}(\bar{x}))$
	\item[] {\bf Compute} BFGS or Newton update $\delta \bar{u}^{(k)}$
	\item[] {\bf Update} $\bar{u}^{(k+1)} = \bar{u}^{(k)} + \delta \bar{u}^{(k)}$ 
	\item[] {\bf Check} for convergence of BFGS or Newton
	\end{itemize}
	\item[] {\bf EndLoop}
	\end{itemize}	
	\caption{FE$^2$ algorithm using Newton's method or BFGS as solver on the macroscopic scale and an NN as surrogate model on the microscopic scale. The complete algorithm is usually embedded in a (dynamic) load stepping scheme for robustness and global convergence.}
	\label{fig:algo2}
\end{figure}

\section{Numerical results}\label{sec:results}

In this section, we aim for a numerical verification that the developed homogenization approach using the beam frame model gives reasonable results and that the machine learning-based surrogate model accurately approximates them in a computationally efficient manner. Additionally, we investigate the convergence behavior of BFGS and Newton's method for different examples. All stress values presented in this section are given in megapascals (MPa).

\subsection{Microscopic simulation and NN prediction}
First, we want to illustrate the behavior of the NN in comparison to the beam frame solver with the homogenization of the Piola-Kirchhoff stresses by presenting some results for microscopic deformations, the resulting Piola-Kirchhoff stress computed with the homogenization method, and the prediction of the NN for the same deformation. We consider the three-dimensional RVE which is also shown in \Cref{fig:micro_results} with the microscopic boundary conditions that are introduced in \Cref{sec:methods}. The results of uniaxial compression tests in each of the three dimensions and two examples which show a shift of the RVE. For the presented microscopic deformations the predictions of the NN are very close to the Piola-Kirchhoff stress tensors computed with the beam frame solver and the homogenization method. \Cref{fig:micro_results} also shows the norm of the deviation between both computed stress tensors. The Frobenius norm is used for the evaluation. 

\begin{table}
	\begin{tabular}{|c|c|c|c|}
		\hline
		deformation gradient & $\overline{P}_{BF}$ & $\overline{P}_{NN}$ & $\frac{\lVert \overline{P}_{BF} - \overline{P}_{NN}\rVert_F}{\lVert \overline{P}_{BF}\rVert_F}$	\\	\hline
		 $\overline{F}=\left(\begin{smallmatrix}
			0.9 & 0 & 0	\\
			0 & 1 & 0 \\
			0 & 0 & 1
		\end{smallmatrix}\right)$& $\left(\begin{smallmatrix}
			-111.53 & 0.53 & 0.43	\\
		0.53 & -2.84 & 1.79 \\
		0.43 & 1.79 & -7.25
	\end{smallmatrix}\right)$ & $\left(\begin{smallmatrix}
	-111.55 & 0.53 & 0.44	\\
	0.54 & -2.83 & 1.76 \\
	0.43 & 1.77 & -7.23
\end{smallmatrix}\right)$ & 3.77e-4\\ \hline 
$\overline{F}=\left(\begin{smallmatrix}
1 & 0 & 0	\\
0 & 0.9 & 0 \\
0 & 0 & 1
\end{smallmatrix}\right)$& $\left(\begin{smallmatrix}
-5.92 & -1.86 & -0.49	\\
-1.86 & -109.60 & -0.12 \\
-0.50 & -0.12 & -6.13
\end{smallmatrix}\right)$ & $\left(\begin{smallmatrix}
-5.93 & -1.86 & -0.50	\\
-1.86 & -109.59 & -0.14 \\
-0.50 & -0.13 & -6.11
\end{smallmatrix}\right)$ & 3.25e-4 \\	\hline
$\overline{F}=\left(\begin{smallmatrix}
1 & 0 & 0	\\
0 & 1 & 0 \\
0 & 0 & 0.9
\end{smallmatrix}\right)$& $\left(\begin{smallmatrix}
-6.98 & 0.74 & -1.45	\\
0.74 & -4.03 & -0.04 \\
-1.45 & -0.04 & -109.46
\end{smallmatrix}\right)$ & $\left(\begin{smallmatrix}
-6.99 & 0.74 & -1.46	\\
0.75 & -4.02 & -0.06 \\
-1.46 & -0.05 & -109.45
\end{smallmatrix}\right)$ & 2.56e-4 \\	\hline
$\overline{F}=\left(\begin{smallmatrix}
1 & 0.2 & 0	\\
0 & 1 & 0 \\
0 & 0 & 1
\end{smallmatrix}\right)$& $\left(\begin{smallmatrix}
1.09 & 104.59 & 0.82	\\
104.59 & 0.15 & -3.30 \\
0.82 & -3.30 & -1.97
\end{smallmatrix}\right)$ & $\left(\begin{smallmatrix}
1.08 & 104.58 & 0.82	\\
104.59 & 0.16 & -3.32 \\
0.83 & -3.32 & -1.95
\end{smallmatrix}\right)$ & 2.52e-4 \\	\hline
$\overline{F}=\left(\begin{smallmatrix}
1 & 0 & 0.2	\\
0 & 1 & 0 \\
0 & 0 & 1
\end{smallmatrix}\right)$& $\left(\begin{smallmatrix}
-4.19 & 0.14 & 107.62	\\
0.14 & -1.32 & 0.09 \\
107.62 & 0.09 & 2.01
\end{smallmatrix}\right)$ & $\left(\begin{smallmatrix}
-4.19 & 0.14 & 107.63	\\
0.14 & -1.29 & 0.07 \\
107.64 & 0.08 & 2.03
\end{smallmatrix}\right)$ & 3.26e-4 \\	\hline
	\end{tabular}
	\caption{Piola-Kirchhoff stress computed for the three-dimensional RVE with beam frame homogenization method and neural network}
	\label{fig:micro_results}
\end{table}

\subsection{Convergence of the macroscopic solver}
As already mentioned we use (Quasi)-Newton methods~\cite{NoceWrig06} for solving the macroscopic problem. Exact Newton methods require the computation of a tangent which is unknown for the case of the beam RVEs and its computation is expected to be expensive. The BFGS method is a common choice to overcome the problem of computing the tangent since it requires only an initial guess for the tangent and adjusts this approximation in each iteration using only first order derivatives, as, e.g., the gradient. In general, this approach lacks the quadratic order of convergence of Newton's method~\cite{Nocedal1987}, but each iteration can be computed much faster. In the present work, we consider Newton's method with an exact tangent for the NN approach and a Quasi-Newton method with an approximate tangent for the beam frame approach. The numerical approximation in this case is computed with central differences. In the following, we will refer to this approach as Newton's method with beam frame microstructure. We also consider the BFGS method~\cite{NoceWrig06,YuHong2002} for both NN and beam frame approach with the Wolfe~\cite{Wolfe1969} condition to control the step length. We initialize the Hessian of the first BFGS step with the exact Hessian in the case of the NN and the approximated Hessian in the beam case. This increases the computation time for the first BFGS step but reduces the total number of iterations for the solver. The stopping criteria for both, the Newton and BFGS are both based on the reduction of the relative residual. If the ratio between the residual of the current update and the residual of the initial value is smaller than $1e-10$ the criterium is satisfied and the solver terminates.

For our first two-dimensional test case we consider a square domain on the macroscopic scale. Dirichlet boundary conditions are applied to the left and right boundaries of the domain and these nodes are shifted in $x$ direction depending on their $y$ value; see~\Cref{fig:bc5_sol} for the macroscopic deformation. The macroscopic problem is discretized using bilinear brick finite elements or linear or quadratic triangular ones. We will use this test case to compare different aspects of our macroscopic solvers. Refining the mesh is straight forward, namely by doubling the numbers of finite elements in $x$ and $y$ direction which leads to a quadratic increase in the degrees of freedom (dofs). To find a solution, the solver has to adjust the position of all considered nodes as can be inferred from the solutions in \Cref{fig:bc5_sol}.

Considering the beam frame model for the microscopic problem, both macroscopic solvers are capable for each considered test case to find a solution that reduces the residual until it reaches the stopping criteria. For larger test cases, an adaptive load stepping has proved to increase the stability of both methods. 

In \Cref{fig:macro_iterations} the number of iterations per load step and the total computation time are presented for the beam frame (BF) model as well as for the surrogate NN on the microscopic scale. For the macroscopic discretization P1 elements are considered. The table shows that on average the BFGS method requires only slightly more iterations per load step than the Newton method. However, due to the reduced computational effort for each iteration the total computation time of the BFGS method is significantly lower. For both microscopic solvers the behavior with respect to the number of iterations is identical since the NN is trained to replicate the deformation behavior and the resulting stresses of the beam frame model. Therefore, the macroscopic solver obtains a very similar information in each quadrature point and likewise computes similar iterations.
Let us remark that the computing times are much lower using the surrogate model but we do not see the factor of 100 or more reported before. This is the case since many MATLAB-based computations, especially on the macroscopic scale, are carried out in both cases in addition to the computation of $\overline{P}$. This portion of the computing time can surely be decreased by a hardware-aware and efficient implementation. Especially the macroscopic solver can either be parallelized in case of Newton's method or transferred to a GPU in case of BFGS. An efficient HPC implementation of our approaches exploiting GPUs is planned for the future.

\begin{table}
	\begin{tabular}{|c|c|ccccc|}
		\hline
		&dofs & $50$ & $162$ & $578$ & $2178$ & $8450$	\\
		\hline
		\multirow{4}{*}{BF solver}&Newton iterations & 2 & 2.25 & 3 & 3.25 & 5 \\
		&BFGS iterations & 3.5 & 3.5 & 4.25 & 5.75 & 9.5	\\
		&Newton comp. time [min] & 0.75 & 2.92 & 14.30 & 61.85 & 459.66 \\
		&BFGS comp. time [min] & 0.27 & 1.04 & 4.36 & 23.55 & 150.49	\\
		\hline
		\multirow{4}{*}{NN solver}&Newton iterations & 2 & 2.25 & 3 & 3.25 & 5 \\
		&BFGS iterations & 3.5 & 3.5 & 4.25 & 5.75 & 9.5	\\
		&Newton comp. time [min] & 0.23 & 0.84 & 3.78 & 14.09 & 87.33  \\
		&BFGS comp. time [min] & 0.06 & 0.20 & 0.70 & 3.01 & 26.33	\\
		\hline
	\end{tabular}
	\caption{Number of macroscopic iterations per load step and total computing time for P1 discretization.}
	\label{fig:macro_iterations}
\end{table}

\subsection{Comparison of different macroscopic grids}
The NN has only been trained on simulation results computed with quadrilateral elements and bilinear Ansatz functions on the macroscopic scale. The numerical test cases use also discretizations with triangular elements and linear (P1) or quadratic (P2) Ansatz functions. To confirm that the NN can be used for P1 and P2 discretizations despite it was not trained on data obtained with P1 or P2 simulations, the deviation of the solutions computed with the beam frame model are compared to the solutions computed with the surrogate model.

\begin{figure*}
	\begin{minipage}{0.3\textwidth}
		\includegraphics[height=0.215\textheight]{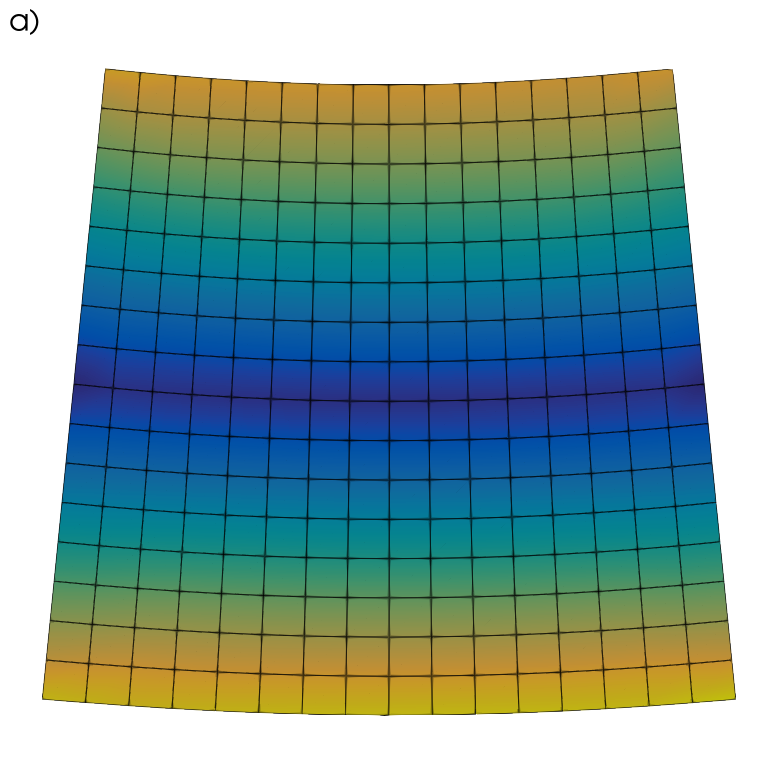}
	\end{minipage}
	\begin{minipage}{0.3\textwidth}
		\includegraphics[height=0.215\textheight]{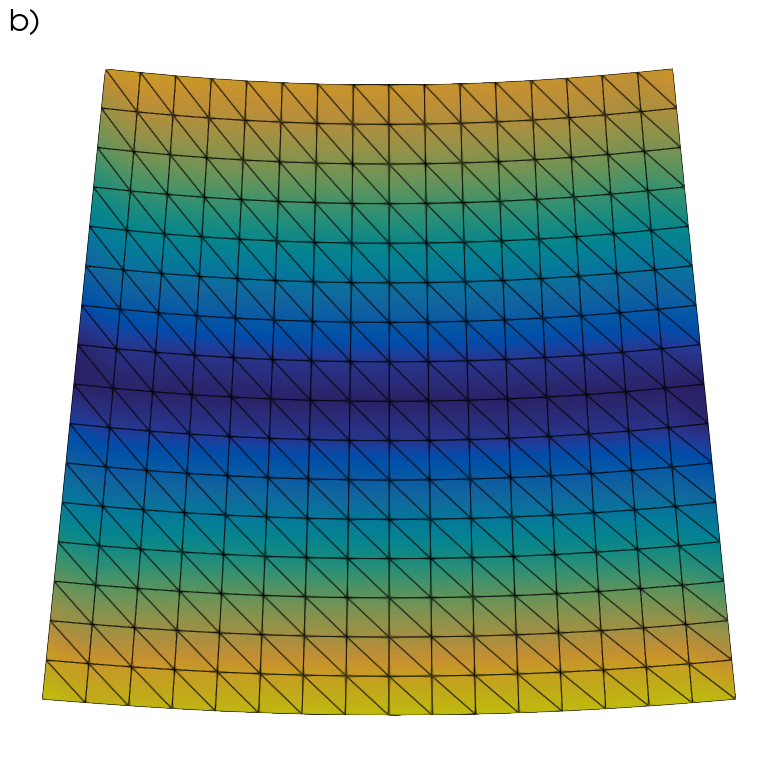}
	\end{minipage}
	\begin{minipage}{0.39\textwidth}
		\includegraphics[height=0.215\textheight]{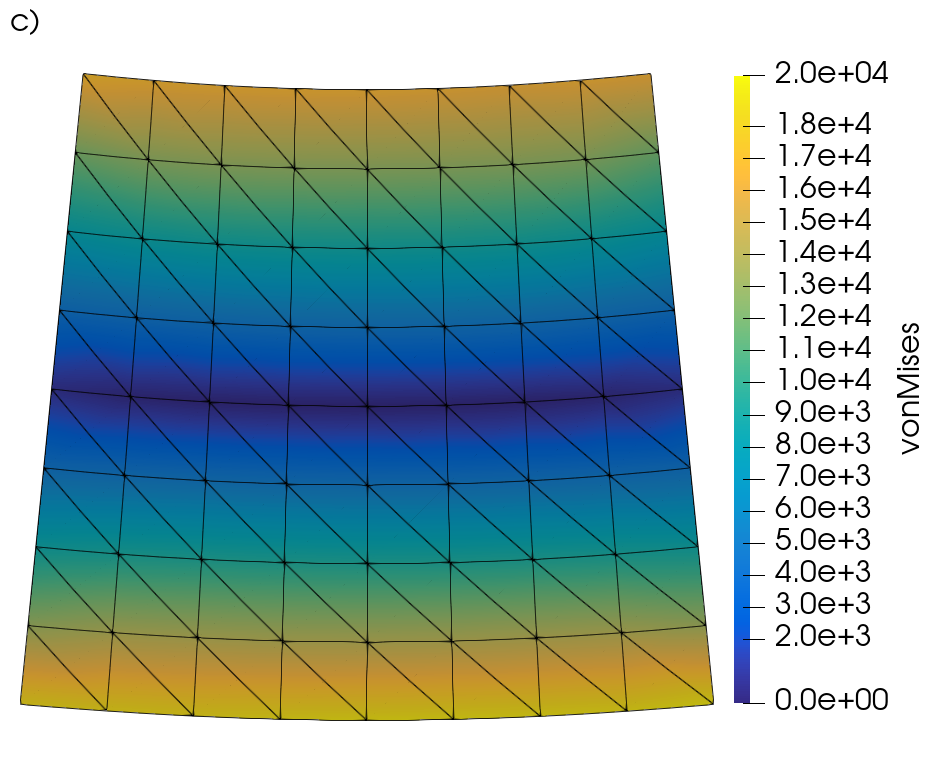}
	\end{minipage}
	\caption{Solutions of the deformation test with Q1 (a), P1 (b), and P2 (c) Ansatz functions.}
	\label{fig:bc5_sol}
\end{figure*}

The deviations of the solutions computed with the NN from the solutions computed with the beam frame model are shown in \Cref{fig:norms}. In the tables $u_{NN}$ refers to the solution computed with the NN and $u_{BF}$ refers to the solution computed with the beam frame model. The expressions $\bar{\sigma}_{vM,NN}$ and $\bar{\sigma}_{vM,BF}$ refer to the corresponding von Mises stresses. We consider different macroscopic meshes with the level of refinement marked by the corresponding degrees of freedom (dofs).

\begin{figure}
	\begin{center}
		\includegraphics[width=0.5\textwidth]{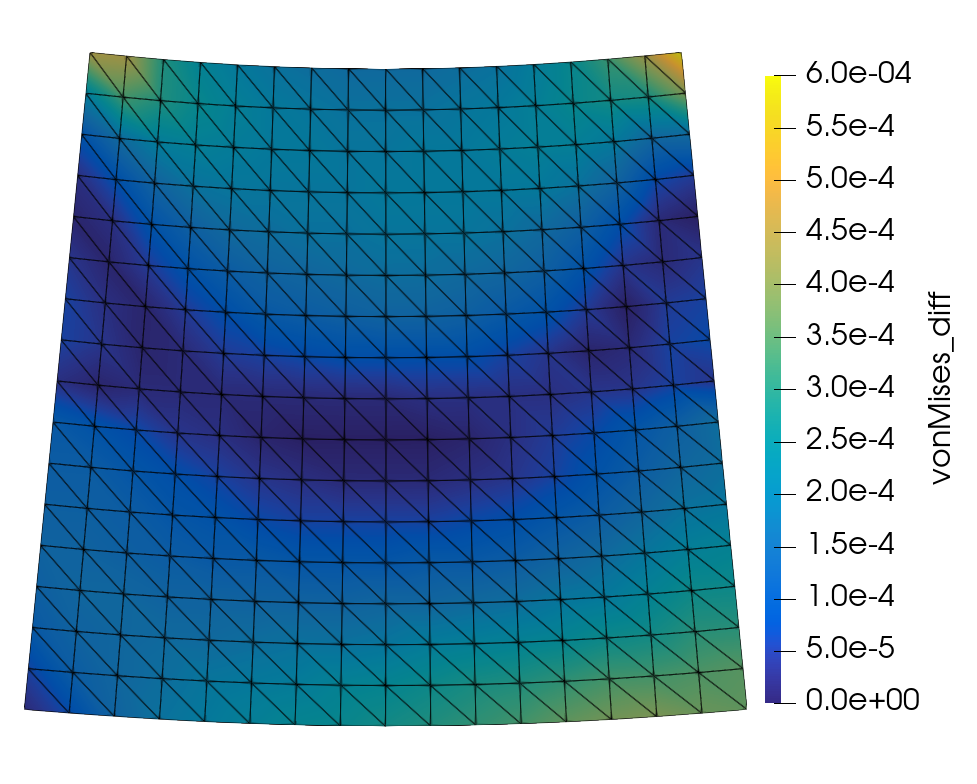}
	\end{center}
	\caption{Difference in Von-Mises stress between the NN solution and the solution computed with the beam frame model.}
	\label{fig:bc5_diff}
\end{figure}

We can observe that the errors of the NN solution computed with different Ansatz functions at the same level of refinement are very similar. The choice of the Ansatz function does not seem to have an influence on the quality of the macroscopic solution when the surrogate model is used. Let us remind that the network has only been trained on data obtained from solutions computed with bilinear square elements. Nonetheless the NN is capable to yield equally good solutions for other element shapes and Ansatz functions. In general, this is expected, since the NN is trained on localized data which does not really depend on the finite element discretization used on the macroscopic scale. Finally, the computed solutions for a refinement with 578 dofs are presented in \Cref{fig:bc5_sol} and the absolute difference between the Von-Mises stress of the solutions computed with the NN and with the beam frame model is shown in \Cref{fig:bc5_diff}.
The macroscopic error of the NN solution is not evenly distributed and higher error values can be observed near the boundaries of the domain. However, the expected distribution of the derivation values is difficult to predict since the error depends locally on the quality of the networks prediction for the Piola-Kirchhoff stress.

The results show that the deviation of the solution computed with the NN from the solution computed with the beam frame model slightly increases with an increasing number of elements. The behavior can be explained with the fact that the computation of the Piola Kirchhoff stress in a single integration point already slightly deviates from the averaged stress computed with the beam frame model. This error is small for each evaluation but in finer macroscopic meshes with more evaluation points these small errors can sum up to a larger macroscopic deviation in the deformation. In the considered mesh refinements the error stays relatively small and for some refinements and norm measurements the error seems to stay constant or even decreases. This behavior leaves room for the assumption that the error might converge with a further refinement and might not increase significantly when a certain level is reached. Unfortunately we are not able to support this assumption with further data at this point. We plan larger simulations using a more efficient software framework than MATLAB in the future.

\begin{table*}
		\begin{tabular}{|lc|ccccc|}
		\hline
		\multicolumn{1}{|c}{} & \multicolumn{1}{c|}{} & \multicolumn{5}{c|}{\bf degrees of freedom}\\
			\multicolumn{1}{|l}{\bf norm} & \multicolumn{1}{c|}{\bf discr.} &  $50$ & $162$ & $578$ & $2178$ & $8450$	\\
			\hline
			& Q1 & $2.09 e-4$ & $2.44 e-4$ & $2.58e-4$ & $2.63e-4$ & $2.65e-4$	\\
			$\frac{\lVert \bar{\mathbf{u}}_{NN} - \bar{\mathbf{u}}_{BF}\rVert_2}{\lVert\bar{\mathbf{u}}_{BF}\rVert_2}$&P1& $2.15 e-4$ & $2.49 e-4$ & $2.63e-4$ & $2.68e-4$ & $2.70e-4$	\\
			&P2& $2.23 e-4$ & $2.52 e-4$ & $2.61e-4$ & $2.64e-4$ & $2.65e-4$	\\\hline
			&Q1 & $4.06 e-4$ & $4.17 e-4$ & $4.23e-4$ & $4.24e-4$ & $4.24e-4$	\\
			$\frac{\lVert \bar{\mathbf{u}}_{NN} - \bar{\mathbf{u}}_{BF}\rVert_\infty}{\lVert\bar{\mathbf{u}}_{BF}\rVert_\infty}$ &P1& $3.70 e-4$ & $4.01 e-4$ & $4.18e-4$ & $4.22e-4$ & $4.24e-4$\\
			&P2& $2.56 e-4$ & $2.79 e-4$ & $2.75e-4$ & $2.74e-4$ & $2.73e-4$	\\\hline\hline
			&Q1& $3.04 e-4$ & $3.09 e-4$ & $3.14e-4$ & $3.14e-4$ & $3.12e-4$	\\
			$\frac{\lVert \bar{\sigma}_{vM,NN} - \bar{\sigma}_{vM,BF}\rVert_2}{\lVert\bar{\sigma}_{vM,BF}\rVert_2}$ &P1& $2.59 e-4$ & $2.85 e-4$ & $3.01e-4$ & $3.08e-4$ & $3.09e-4$\\
			&P2& $3.06 e-4$ & $3.52 e-4$ & $3.69e-4$ & $3.73e-4$ & $3.73e-4$	\\\hline
			&Q1& $4.23 e-4$ & $4.61 e-4$ & $7.10e-4$ & $1.01e-3$ & $1.34e-3$	\\
			$\frac{\lVert \bar{\sigma}_{vM,NN} - \bar{\sigma}_{vM,BF}\rVert_\infty}{\lVert\bar{\sigma}_{vM,BF}\rVert_\infty}$ &P1& $2.94 e-4$ & $3.62 e-4$ & $6.04e-4$ & $9.04e-4$ & $1.24e-3$\\
			&P2& $3.68 e-4$ & $4.97 e-4$ & $8.27e-4$ & $1.14e-3$ & $1.47e-3$	\\\hline
			
		\end{tabular}
	\caption{Deviations of the NN solutions and the BF solutions in different norms for different macroscopic discretizations (discr.), that is, bilinear brick (Q1) and linear (P1) and quadratic (P2) triangular elements.}
		\label{fig:norms}
\end{table*}

\subsection{Compression of a punched plate geometry}
As a more complex example, we consider a two-dimensional plate with a circular hole in the middle. The square shaped geometry has an edge length of 1 and the hole center has a diameter of 0.4. We want to apply a compression test on this geometry to compare the macroscopic solutions using the beam frame approach and the NN. The geometry is discretized with P1 elements for both microscopic solvers and the macroscopic solution is computed with the BFGS method. For the compression test, the left and right boundaries are shifted towards the hole and fixed in $x$ and $y$ direction with Dirichlet boundary conditions.

\begin{figure*}[t]
	\begin{minipage}{0.49\textwidth}
		\includegraphics[width=\textwidth]{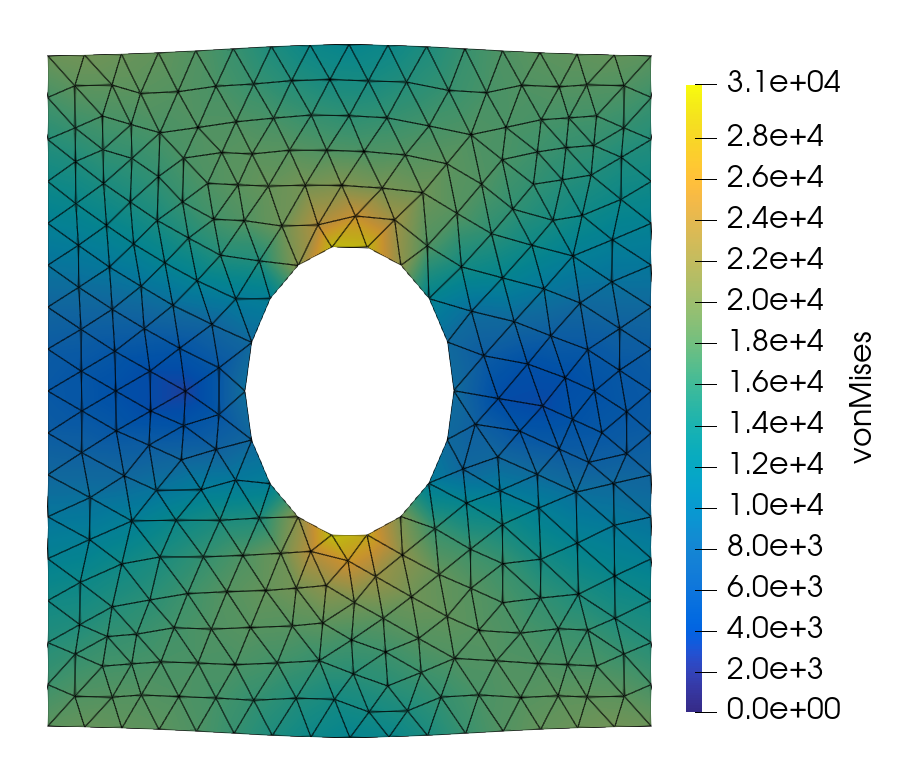}
	\end{minipage}
	\begin{minipage}{0.49\textwidth}
	\includegraphics[width=\textwidth]{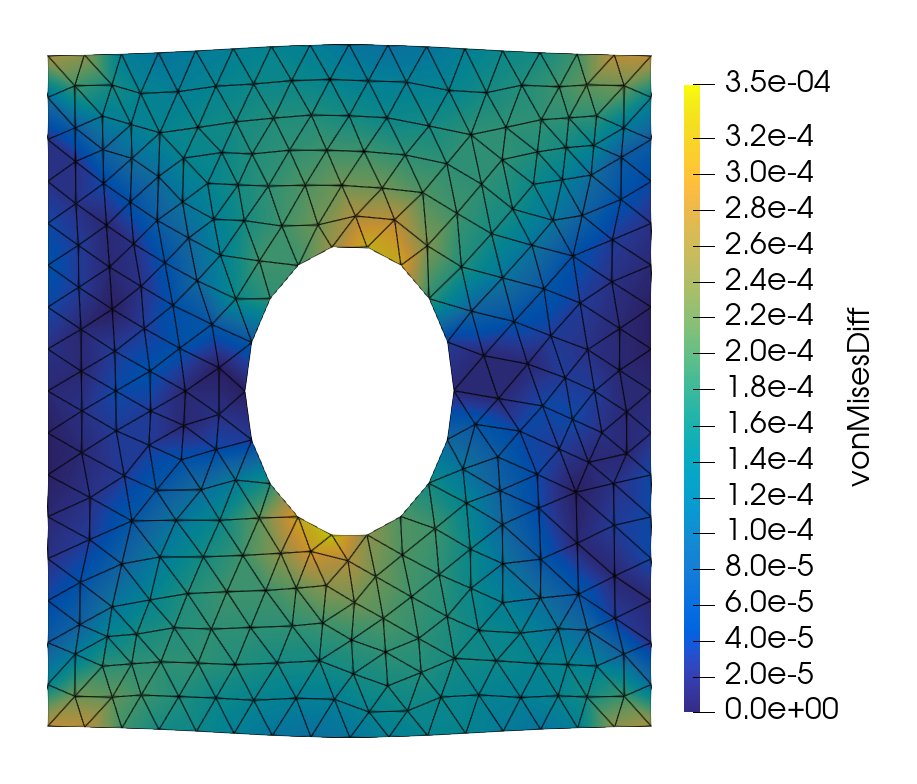}
	\end{minipage}
	\caption{Compressed punched plate computed with the NN approach and P1 elements. Coloring of the elements in the left image represents the resulting Von-Mises stresses and coloring in the right image shows the absolute difference of the Von-Mises stress to the solution computed with the beam frame solver.}
	\label{fig:hole2d}
\end{figure*}

The result of the compression test computed with the NN approach and the resulting difference in Von-Mises stress are shown in \Cref{fig:hole2d}. The largest stress differences between the two solutions are observed in the regions of the geometry with also the highest stress values. The total error between the macroscopic solutions is evaluated with the following norms:

\begin{align*}
	\frac{\lVert \bar{\mathbf{u}}_{NN} - \bar{\mathbf{u}}_{BF}\rVert_2}{\lVert\bar{\mathbf{u}}_{BF}\rVert_2} &= 1.81e-4,	\\
	\frac{\lVert \bar{\mathbf{u}}_{NN} - \bar{\mathbf{u}}_{BF}\rVert_\infty}{\lVert\bar{\mathbf{u}}_{BF}\rVert_\infty} &= 2.30e-4,	\\
	\frac{\lVert \bar{\sigma}_{vM,NN} - \bar{\sigma}_{vM,BF}\rVert_2}{\lVert\bar{\sigma}_{vM,BF}\rVert_2} &= 2.86e-4,	\\
	\frac{\lVert \bar{\sigma}_{vM,NN} - \bar{\sigma}_{vM,BF}\rVert_\infty}{\lVert\bar{\sigma}_{vM,BF}\rVert_\infty} &= 3.47e-4.
\end{align*}

\subsection{Compression of a punched plate geometry with five holes}
\begin{figure*}[b]
	\begin{minipage}{0.495\textwidth}
		\includegraphics[width=\textwidth]{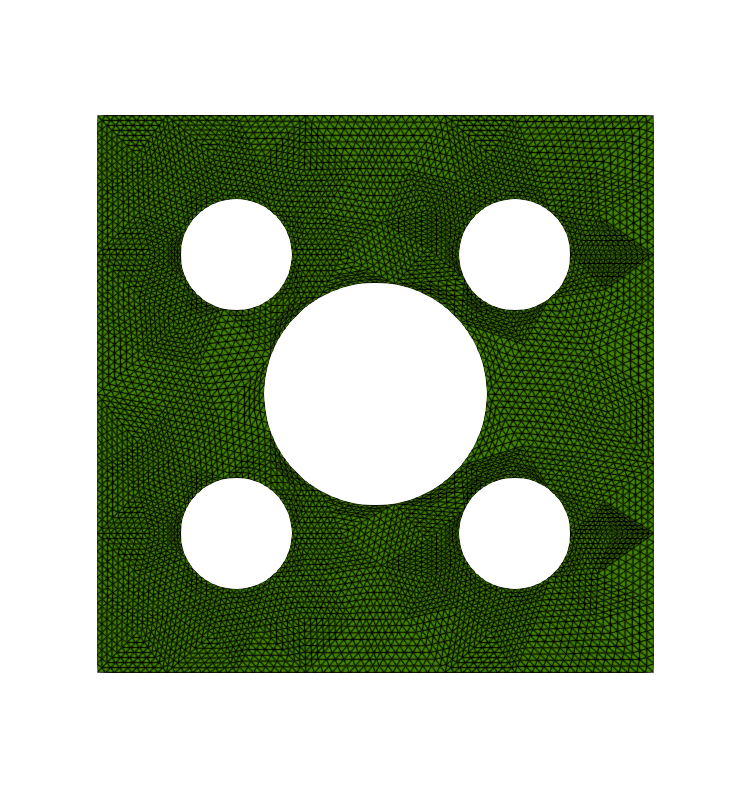}
	\end{minipage}
	\begin{minipage}{0.5\textwidth}
	\includegraphics[width=\textwidth]{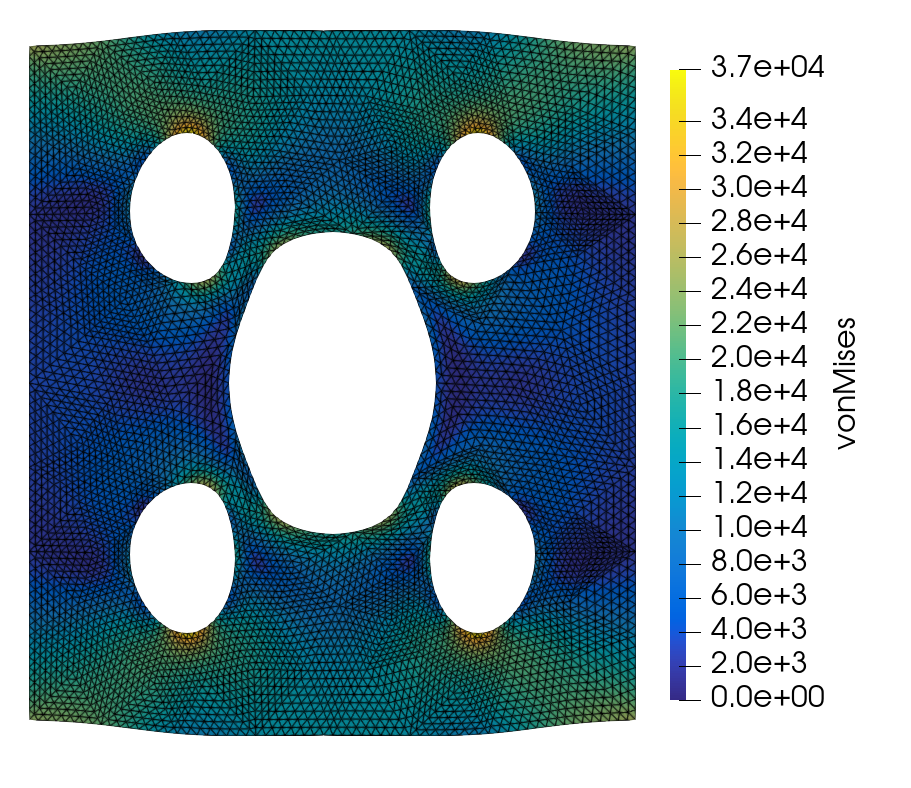}
	\end{minipage}
	\caption{Compressed plate with five punched holes computed with the NN approach and P1 elements. Colors of the elements represent the resulting Von-Mises stress.}
	\label{fig:holes2d}
\end{figure*}

The NN approach allows the simulation of finer and more complex geometries without exhausting the computational resources. As an example we consider a two-dimensional square-shaped plate with an edge length of 1 and five circular holes. One hole is in the center of the plate with a diameter of 0.4 and four smaller holes with a diameter of 0.2 are placed around it with the same distance from the center. The circles which form the holes are placed without any overlap. For this geometry we apply the same two-dimensional compression test as in the previous example with the left and right boundaries being fixed with Dirichlet boundary conditions. We use linear (P1) elements for a fine discretization of the geometry which yields 23,824 degrees of freedom. The computed BFGS solution of this deformation is presented in \Cref{fig:holes2d} with the resulting Von-Mises stresses expressed by the coloring of the elements.

\subsection{Comparison of three-dimensional test cases}

\begin{figure*}
	\begin{minipage}{0.32\textwidth}
		\includegraphics[width=\textwidth]{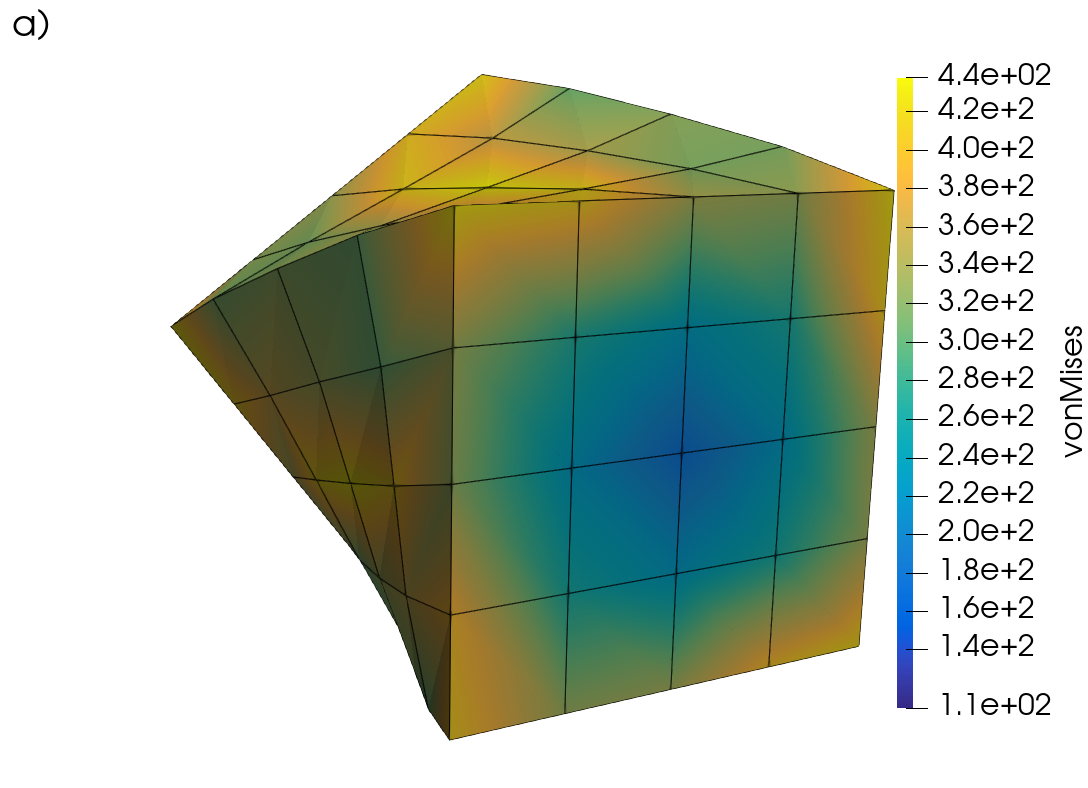}
	\end{minipage}
	\begin{minipage}{0.32\textwidth}
		\includegraphics[width=\textwidth]{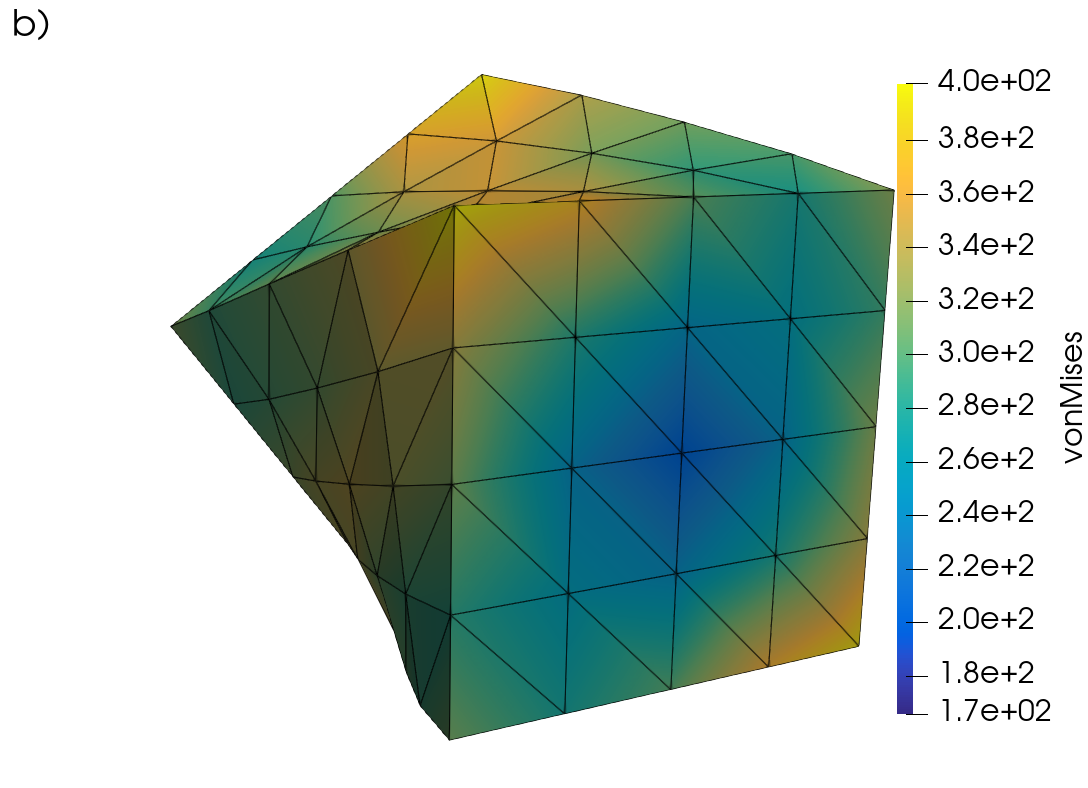}
	\end{minipage}
	\begin{minipage}{0.32\textwidth}
		\includegraphics[width=\textwidth]{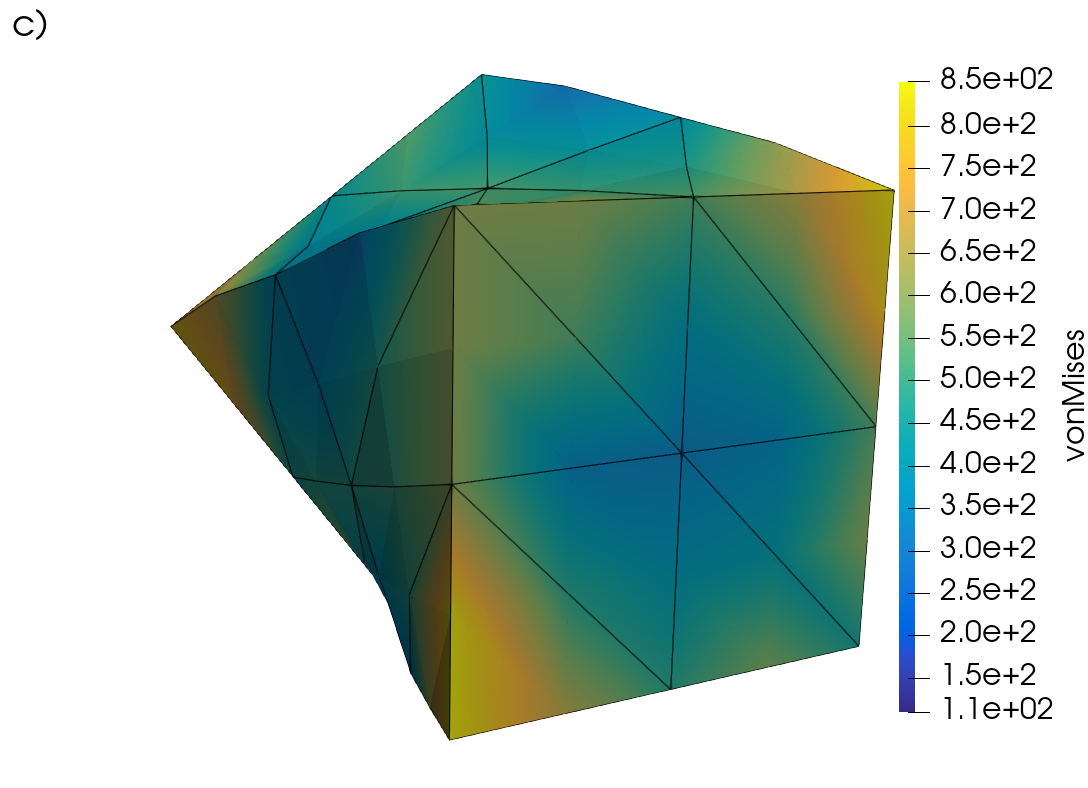}
	\end{minipage}
	\caption{Solution of the three-dimensional deformation test computed with Q1 (a), P1 (b), and P2 (c) Ansatz functions. Coloring of the elements represents the resulting Von-Mises stresses.}
	\label{fig:sol3d}
\end{figure*}

For our macroscopic test problem in three dimensions we consider a regular cube. Dirichlet boundary conditions are applied to the boundaries in $x$ direction and the affected nodes are both rotated around the center of the face and pulled in $x$ direction away from the center of the cube. The cube is rotated by a total of 36 degrees. This geometry again fits well for a comparison of different finite element grids since trilinear cube elements can be used for the discretization as well as linear or quadratic tetrahedral elements. The refinement of the mesh is achieved by doubling the number of elements in each dimension and yield a cubic increase in the dofs. Solutions of the test problem computed with the different Ansatz functions are presented in \Cref{fig:sol3d}. A comparison of the error between the macroscopic solution computed with the beam frame model and the solution resulting from the surrogate approach is shown in \Cref{fig:norms3d}. Multiple levels of refinement are computed and the results show that the norms of the deviations are in a very similar range to the norms of the regarded two-dimensional test case. The table shows that the errors for one refinement level is about the same for all Ansatz functions and does not indicate a better performance of the NN approach for one of the discretizations. For higher numbers of degrees of freedom the errors increase slightly, without occurrence of any significant outliers.

\begin{table}
		
	\begin{tabular}{|lc|ccc|}
		\hline
		\multicolumn{1}{|c}{} & \multicolumn{1}{c|}{} & \multicolumn{3}{c|}{\bf degrees of freedom}\\
		\multicolumn{1}{|l}{\bf norm} & \multicolumn{1}{c|}{\bf discr.} &  $108$ & $525$ & $3159$ 	\\
		\hline
		& Q1 & $6.32 e-5$ & $8.99 e-5$ & $1.09e-4$	\\
		$\frac{\lVert \bar{\mathbf{u}}_{NN} - \bar{\mathbf{u}}_{BF}\rVert_2}{\lVert\bar{\mathbf{u}}_{BF}\rVert_2}$&P1& $5.66 e-5$ & $8.68 e-5$ & $1.07e-4$	\\
		&P2& $3.38 e-5$ & $8.50 e-5$ & $1.06e-4$	\\\hline
		&Q1 & $1.78 e-4$ & $1.79 e-4$ & $1.52e-4$	\\
		$\frac{\lVert \bar{\mathbf{u}}_{NN} - \bar{\mathbf{u}}_{BF}\rVert_\infty}{\lVert\bar{\mathbf{u}}_{BF}\rVert_\infty}$ &P1& $1.18 e-4$ & $1.52 e-4$ & $1.42e-4$	\\
		&P2& $7.28 e-5$ & $1.32 e-4$ & $1.41e-4$	\\\hline\hline
		&Q1& $1.29 e-4$ & $8.97 e-5$ & $8.63e-5$	\\
		$\frac{\lVert \bar{\sigma}_{vM,NN} - \bar{\sigma}_{vM,BF}\rVert_2}{\lVert\bar{\sigma}_{vM,BF}\rVert_2}$ &P1& $6.23 e-5$ & $5.81 e-5$ & $6.38e-5$	\\
		&P2& $1.15 e-4$ & $1.54 e-4$ & $2.96e-4$	\\\hline
		&Q1& $3.55 e-4$ & $2.43 e-4$ & $2.11e-4$\\
		$\frac{\lVert \bar{\sigma}_{vM,NN} - \bar{\sigma}_{vM,BF}\rVert_\infty}{\lVert\bar{\sigma}_{vM,BF}\rVert_\infty}$ &P1& $1.18 e-4$ & $1.17 e-4$ & $1.49e-4$	\\
		&P2& $2.28 e-4$ & $3.63 e-4$ & $1.37e-3$\\\hline
		
	\end{tabular}
	\caption{Deviations of the NN solutions and the BF solutions in different norms for different macroscopic discretizations (discr.), that is, bilinear brick (Q1) and linear (P1) and quadratic (P2) triangular elements.}
	\label{fig:norms3d}
\end{table}

\subsection{Torsion of a cube with a cylindrical hole}

\begin{figure}
	\begin{minipage}{0.49\textwidth}
		\includegraphics[width = \textwidth]{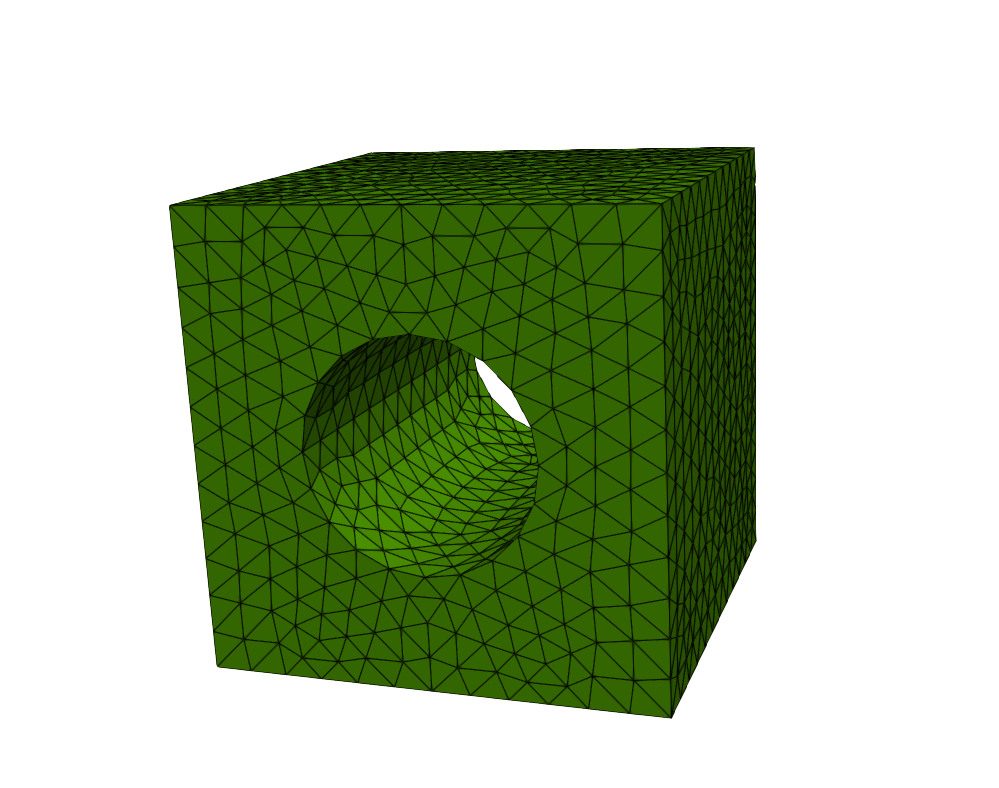}
	\end{minipage}
\begin{minipage}{0.49\textwidth}
	\includegraphics[width = \textwidth]{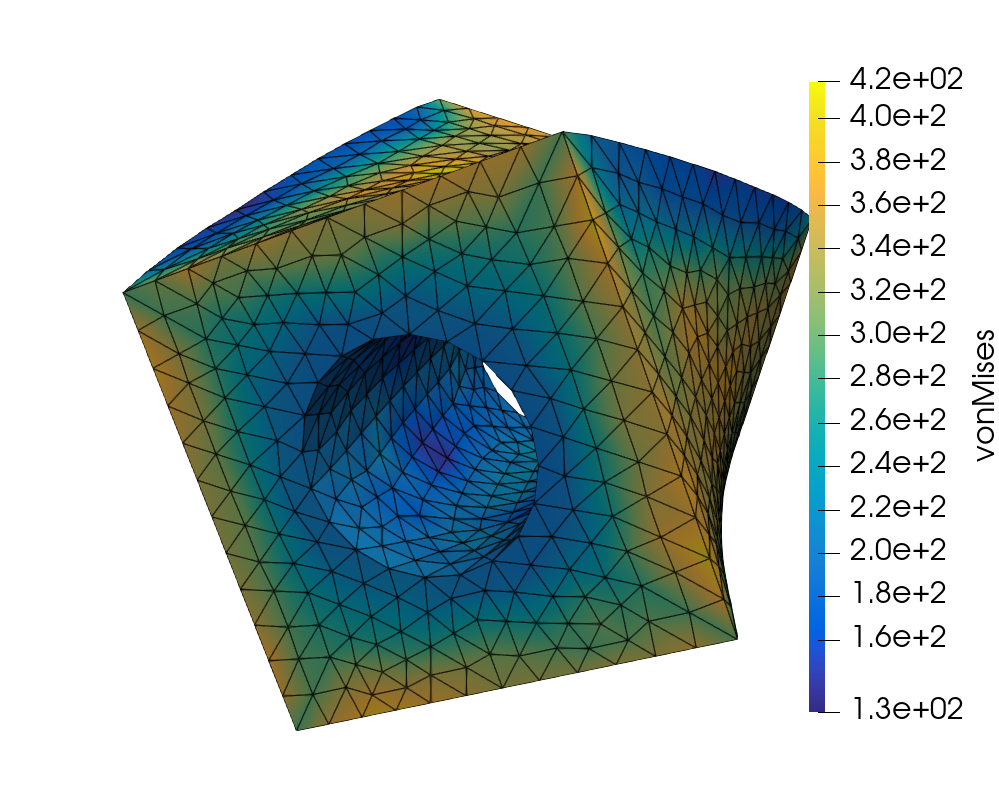}
\end{minipage}
\caption{Torsion of a cube geometry with a circular hole discretized with P1 elements.}
\label{fig:hole3d_torsion}
\end{figure}

As a more complex example of a three-dimensional deformation, we consider a cube with an edge length of 1 and a cylindrical hole in the middle. The cylinder-shape is placed parallel to the $x$ axis with a diameter of 0.5 and the geometry is discretized with linear P1 elements. Dirichlet boundary conditions are applied to the boundaries with $x$ is equal to zero or one and similar to the previous test example each affected node is twisted and shifted in $x$ direction away from the cubes center. We again consider a total rotation of 36 degrees. The result with the corresponding Von-Mises stresses are presented in \Cref{fig:hole3d_torsion}.

\subsection{Torsion of a cylinder}

\begin{figure*}
	\begin{minipage}{0.49\textwidth}
		\includegraphics[width = \textwidth]{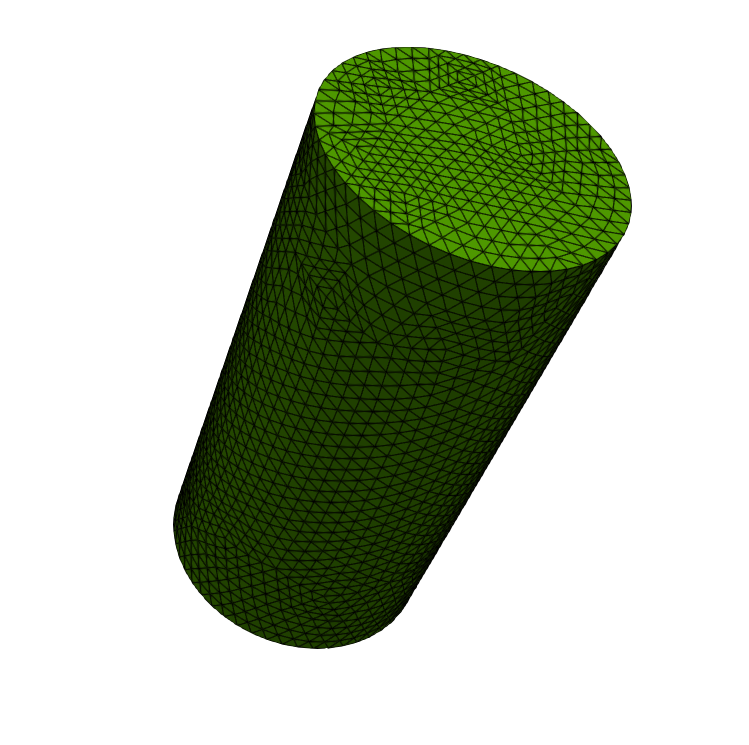}
	\end{minipage}
\begin{minipage}{0.49\textwidth}
\includegraphics[width = \textwidth]{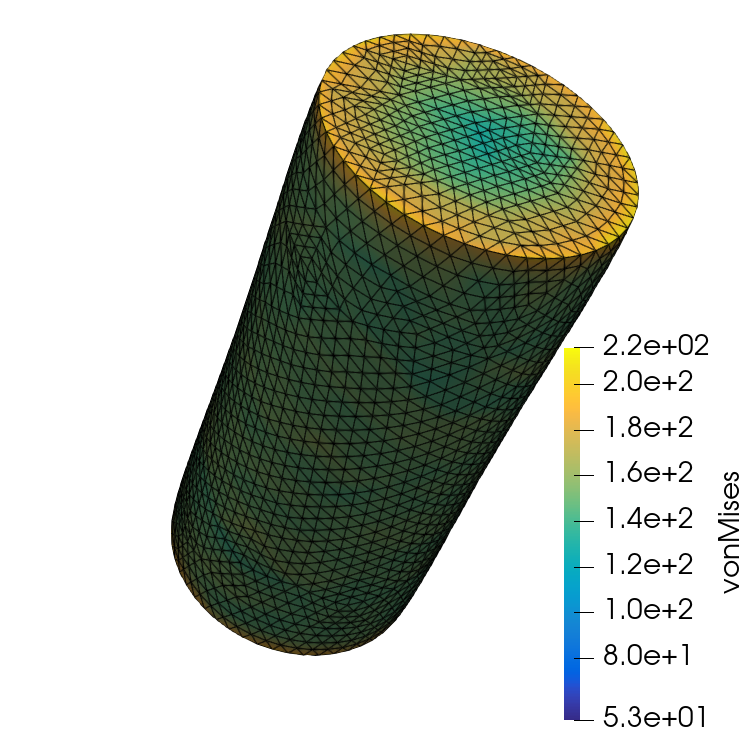}
\end{minipage}
	\caption{Torsion of a cylinder geometry discretized with P1 elements. Left image shows the undeformed geometry and right image shows the computed BFGS solution.}
	\label{fig:cyl3d_torsion}
\end{figure*}

As another example of a three-dimensional deformation we consider a cylinder with the length of 2 and a diameter of 1. The cylinder-shape is placed parallel to the $x$ axis and the geometry is discretized with linear P1 elements. The resulting discretization has 27,231 dofs. For the torsion test case Dirichlet boundary conditions are applied to the boundaries with $x$ equals zero or one similar to the previous test example each affected node is twisted and shifted in $x$ direction away from the cylinders center. The total rotation of the geometry is 36 degrees like in the previous test cases. The results in \Cref{fig:cyl3d_torsion} show that the distribution of the Von-Mises stresses is more evenly compared to the examples of the torsion test of a cube in \Cref{fig:sol3d}.

\bmhead{Acknowledgments}

The authors gratefully acknowledge the Program Directorate Transport of the German Aerospace Centre for funding the Road Traffic Project FFAE. We would also like to thank Shivangi Aney (German Aerospace Centre) for providing the software to generate the open-porous RVEs and Simon Klaes (University of Cologne) for supporting the visualization of microscopic and macroscopic solutions in Paraview.

\begin{appendices}
	\section{Derivation of determinate and inverse for the beam frame model}\label{appendix1}
	Based on the assumptions of the beam frame model as described in section \Cref{sec:methods} the deformation of each beam element can be expressed with respect to the distance of its starting node.
	
	\begin{align*}
		\mathbf{u}(x,y,z) &= \tilde{\mathbf{u}}(\xi) = \left(\begin{array}{c}
			a_x\xi^3+b_x\xi^2+c_x\xi+d_x  \\
			a_y\xi^3+b_y\xi^2+c_y\xi+d_y  \\
			a_z\xi^3+b_z\xi^2+c_z\xi+d_z
		\end{array}\right)
	\end{align*}
	with  $\xi = \frac{\lVert (x,y,z)^T-v_1\rVert}{\lVert v_2-v_1\rVert}\in[0,1]$ being the relative distance from the coordinate $X$ to the first vertex of the beam $v_1$. With this expression it is possible to compute the Jacobian of $u$ as
	\begin{equation*}
		\nabla \mathbf{u}(X) = \left(\begin{array}{ccc}
			u_{xx} & u_{yx} & u_{zx}  \\
			u_{xy} & u_{yy} & u_{zy}  \\
			u_{xz} & u_{yz} & u_{zz}
		\end{array}\right)
	\end{equation*}
	with $u_{x_i x_j} = \frac{\partial \tilde{u}_{x_i}}{\partial \xi} \frac{\partial \xi}{\partial x_j}= \left(3a_{x_i}\xi^2+2b_{x_i}\xi+c_{x_i}\right) \frac{(v_2-v_1)_j}{\lVert v_2-v_1\rVert^2}$. This elementwise notation is equivalent to the dyadic product $\nabla \mathbf{u}(X) = \nabla_\xi\, \tilde{\mathbf{u}}(\xi) \otimes \nabla \xi$. Therefore, the deformation gradient can be expressed in the form 
	\begin{equation*}
		F = I + p\otimes q.
	\end{equation*}
	for $p = \nabla_\xi\, \hat{u}(\xi)$ and $q=\nabla \xi$. 
	
	For the three-dimensional case we will now derive that for every matrix $F\in\mathbb{R}^{3\times 3}$ that follows this relation and any given $p,q\in\mathbb{R}^3$ the determinant is given by $\det(F)=1+\text{trace}(p\otimes q)$. Using Leibniz formula for the computation of the determinant of the $3\times 3$ matrix yields
	\begin{align*}
		\det(F) =& (1+p_x\,q_x)\cdot(1+p_y\,q_y)\cdot(1+p_z\,q_z) + 2\,p_x\,q_x\,p_y\,q_y\,p_z\,q_z	\\
		&-(1+p_x\,q_x)\cdot p_y\,q_y\,p_z\,q_z-(1+p_y\,q_y)\cdot p_x\,q_x\,p_z\,q_z	\\
		&-(1+p_z\,q_z)\cdot p_x\,q_x\,p_y\,q_y.
	\end{align*}
	After resolving the brackets most of the terms cancel out which leads to the simplified form
	\begin{equation*}
		\det(F)=1+p_x\,q_x\cdot p_y\,q_y\cdot p_z\,q_z = 1+\text{trace}(p\otimes q).
	\end{equation*}
	With the given expression for the determinant it is also possible to derive the formulation for the inverse $F^{-1} = \frac{1}{\det(F)}\left((1+\det(F))-F\right)$. For the computation of the inverse we use the well-known relation of the inverse with the adjugate matrix $F^{-1}=\frac{1}{\det(F)}\cdot$ adj$(F)$ with the adjugate matrix given by
	\begin{equation*}
		\text{adj}(F) = \left(\begin{array}{ccc}
			1+p_y\,q_y+p_z\,q_z & -p_x\,q_y & -p_x\,q_z	\\
			-p_y\,q_x & 1+p_x\,q_x+p_z\,q_z & -p_y\,q_z	\\
			-p_z\,q_x & -p_z\,q_y & 1+p_x\,q_x+p_y\,q_y
		\end{array}\right).
	\end{equation*}
	Finally, together with the relation $\det(F)=1+p_x\,q_x\cdot p_y\,q_y\cdot p_z\,q_z$ it is possible to express the inverse matrix as
	\begin{align*}
		F^{-1} = \frac{1}{\det(F)}\text{adj}(F) = \frac{1}{\det(F)}\left((1+\det(F))\cdot I-F\right).
	\end{align*}
	
	\section{Distribution of the training data for the three-dimensional NN}\label{appendix2}
	\begin{table*}[h!]
		\begin{tabular}{|c|ccccccccc|}
			\hline
			variable name & $\overline{F}_{xx}$ & $\overline{F}_{yx}$& $\overline{F}_{zx}$ & $\overline{F}_{xy}$ & $\overline{F}_{yy}$ & $\overline{F}_{zy}$ & $\overline{F}_{xz}$ & $\overline{F}_{yz}$ & $\overline{F}_{zz}$	\\
			\hline
			mean & 0 & -0.03 & -0.01 & -0.02 & -0.02 & 0.01 & 0.02 & -0.01 & 0  \\
			standard deviation & 0.14 & 0.14 & 0.14 & 0.14 & 0.14 & 0.14 & 0.14 & 0.14 & 0.14 \\
			minimum value & -0.91 & -1.53 & -1.29 & -1.12 & -1.52 & -1.39 & -1.09 & -1.48 & -1.35	\\
			maximum value & 1.13 & 1.02 & 1.17 & 1.21 & 0.95 & 1.32 & 1.19 & 1.05 & 1.29	\\
			\hline
		\end{tabular}
		\caption{Description of the input data for the training of the three-dimensional neural network.}
		\label{fig:input_data_tab3D}
	\end{table*}
	
	\begin{table*}[h!]
		\begin{tabular}{|c|ccccccccc|}
			\hline
			variable name & $\overline{P}_{xx}$ & $\overline{P}_{yx}$& $\overline{P}_{zx}$ & $\overline{P}_{xy}$ & $\overline{P}_{yy}$ & $\overline{P}_{zy}$ & $\overline{P}_{xz}$ & $\overline{P}_{yz}$ & $\overline{P}_{zz}$	\\
			\hline
			mean & 0.01 & -0.03 & 0.01 & -0.03 & -0.02 & 0 & 0.01 & 0 & 0  \\
			standard deviation & 0.16 & 0.11 & 0.11 & 0.11 & 0.15 & 0.10 & 0.11 & 0.10 & 0.16 \\
			minimum value & -1.11 & -1.39 & -1.13 & -1.39 & -1.20 & -1.25 & -1.22 & -1.24 & -1.16	\\
			maximum value & 1.81 & 1.20 & 1.49 & 1.20 & 1.75 & 1.35 & 1.49 & 1.36 & 2.00	\\
			\hline
		\end{tabular}
		\caption{Description of the output data for the training of the three-dimensional neural network.}
		\label{fig:output_data_tab3D}
	\end{table*}
\end{appendices}

\bibliography{sn-bibliography}

\end{document}